\documentclass{amsart}
\usepackage{amsmath,amssymb,amscd,latexsym,enumitem,amsthm, hyperref, bbding, fixme}

\usepackage[all]{xy}
\usepackage[T1]{fontenc}

\newcommand{\C}{{\mathbb{C}}}
\newcommand{\N}{{\mathbb{N}}}
\newcommand{\R}{{\mathbb{R}}}
\newcommand{\Z}{{\mathbb{Z}}}

\newcommand{\dimrok}{\mathrm{dim}_{\mathrm{Rok}}}

\newcommand{\halg}{ {\mathrm C}^*(H)}
\newcommand{\salg}{ \mathrm{C}^*(S)}
\newcommand{\ualg}{ \mathrm{C}^*(U)}

\DeclareMathOperator{\G}{\mathcal{G}}
\DeclareMathOperator{\id}{\textup{id}}
\DeclareMathOperator{\into}{\hookrightarrow}

\DeclareMathOperator{\uhf}{\mathcal{U}}

\DeclareMathOperator{\Aut}{Aut}
\DeclareMathOperator{\homeo}{Homeo}

\DeclareMathOperator{\tsum}{\textstyle{\sum}}

\newcounter{number}[section]

\newenvironment{nummer}{\refstepcounter{number}{\bf \noindent\arabic{section}.\arabic{number}}}{}

\newcommand{\bn}{\noindent \begin{nummer} \rm}
\newcommand{\en}{\end{nummer}}

\newenvironment{ntheorem}{\noindent {\bf Theorem:} \it}{}
\newenvironment{nlemma}{\noindent {\bf Lemma:} \it}{}
\newenvironment{nprop}{\noindent {\bf Proposition:} \it}{}
\newenvironment{ndefn}{\noindent {\bf Definition:} }{}
\newenvironment{ncor}{\noindent {\bf Corollary:} \it}{}

\newenvironment{nremark}{\noindent {\bf  Remark:}}{}
\newenvironment{nremarks}{\noindent {\bf Remarks:}}{}

\newenvironment{nexample}{\noindent {\bf  Example:} }{}

\newenvironment{nproof}{\noindent {\emph{Proof.}} }{\mbox{}\hfill 
\rule[-.2ex]{.25em}{1.8ex}}

\title[Group actions on Smale space $\mathrm{C}^*$-algebras]{Group actions on Smale space $\mathrm{C}^*$-algebras}
\author{Robin J. Deeley}
\address{Robin J. Deeley,   Department of Mathematics,
University of Colorado Boulder
Campus Box 395,
Boulder, CO 80309-0395, USA }
\email{robin.deeley@gmail.com}
\author{Karen R. Strung}
\address{Karen R. Strung, Institute for Mathematics, Astrophysics, and Particle Physics, Radboud University, Postbus 9010, 6500 GL Nijmegen,
Netherlands}
\email{k.strung@math.ru.nl}
\date{\today}
\thanks{The second listed author was supported by an IMPACT fellowship cofunded by Ministry of Science and Higher Education grant 3038/7.PR/2014/2 and EC grant PCOFUND-GA-2012-600415, and the Sonata 9 NCN grant 2015/17/D/ST1/02529}
\begin{document}

\begin{abstract}
Group actions on a Smale space and the actions induced on the $\mathrm{C}^*$-algebras associated to such a dynamical system are studied. We show that an effective action of a discrete group on a mixing Smale space produces a strongly outer action on the homoclinic algebra. We then show that for irreducible Smale spaces, the property of finite Rokhlin dimension passes from the induced action on the homoclinic algbera to the induced actions on the stable and unstable $\mathrm{C}^*$-algebras. In each of these cases, we discuss the preservation of properties---such as finite nuclear dimension, $\mathcal{Z}$-stability, and classification by Elliott invariants---in the resulting crossed products.
\end{abstract}

\maketitle

\section{Introduction}

Topological dynamical systems have long been a source for constructing interesting examples of $\mathrm{C}^*$-algebras. Via the Gelfand transform, compact Hausdorff spaces and continuous maps are in one-to-one contravariant correspondence to unital commutative $\mathrm{C}^*$-algebras and unital $^*$-homomorphisms. To study topological dynamical systems using $\mathrm{C}^*$-algebraic techniques, loosely speaking, one encodes dynamical aspects of the system in a noncommutative $\mathrm{C}^*$-algebra. 

For example, to encode the orbit equivalence classes of a topological dynamical system, one uses the transformation group construction to produce a so-called crossed product $\mathrm{C}^*$-algebra. These have been widely studied for minimal dynamical systems because in such a case, the associated $\mathrm{C}^*$-algebra is simple. When the space is a Cantor set with a single automorphism, it was shown by Putnam that the $\mathrm{C}^*$-algebras are classifiable by $K$-theoretic data \cite{Putnam:MinHomCantor}. Moreover, he proved that they are all approximately circle algebras. Subsequently, Giordano, Putnam and Skau showed that there is an isomorphism of the corresponding $\mathrm{C}^*$-algebras if and only if the systems are strong orbit equivalent \cite{GioPutSkau:orbit}.

The question of classification of such $\mathbb{Z}$-actions proved much more difficult for higher dimensional spaces, but was eventually settled by Lin in \cite{Lin:MinDyn} following many partial results (see for example \cite{LinPhi:MinHom, TomsWinter:PNAS, Win:ClassCrossProd, Str:XxSn}): the $\mathrm{C}^*$-algebras associated to minimal dynamical systems on compact metric spaces with finite covering dimension are classifiable by the Elliott invariant for simple separable unital nuclear $\mathrm{C}^*$-algebras. For a given $\mathrm{C}^*$-algebra, the Elliott invariant consists of its ordered $K$-theory, tracial state space, and a pairing map between these two objects. Unlike in the zero-dimensional case, however, how this classification result relates to the dynamical systems involved is unclear. 

In the present paper, we are interested in nonminimal systems called Smale spaces. These were defined by Ruelle \cite{Rue:ThermForm} based on the behaviour of the basic sets associated to Smale's Axiom A diffeomorphisms \cite{Sma:DiffDynSys}. In the present paper, a dynamical system is a compact metric space with a self-homeomorphism. A Smale space is a dynamical system $(X, \varphi)$ which has a local hyperbolic structure: at every point $x \in X$ there is a small neighbourhood which decomposes into the product of a stable and unstable set.  For a Smale space, the dynamical behaviour which we seek to study is the asymptotic behaviour of points. The appropriate $\mathrm{C}^*$-algebras, defined by Putnam \cite{Put:C*Smale} following earlier work by Ruelle \cite{Rue:NCAlgs} encode this behaviour via  the groupoid $\mathrm{C}^*$-algebras,  $\salg$, $\ualg$ and $\halg$, associated to stable, unstable, and homoclinic equivalence relations, respectively. 

When a Smale space is mixing, these three $\mathrm{C}^*$-algebras are each simple, separable and nuclear. Furthermore, the homoclinic $\mathrm{C}^*$-algebra is unital. From the point of view of the structure and classification of $\mathrm{C}^*$-algebras, the path of investigation has shared many similarities to the case of minimal systems: the first, and most successful, results come from the zero dimensional case. Here, a Smale space is a shift of finite type and the associated $\mathrm{C}^*$-algebras are all approximately finite (AF) algebras, which can be classified by $K$-theory. Moreover, two shifts of finite type are eventually conjugate if and only if the associated $\mathrm{C}^*$-algebras are equivariantly isomorphic with respect to the shift maps \cite[Theorems 7.5.7 and 7.5.15]{LindMarcus:SDandCod}.

 As in the case of minimal systems, classification for homoclinic $\mathrm{C}^*$-algebras of mixing Smale spaces in higher-dimensional settings proved to be more involved. However, the authors of the present paper showed in \cite{DS:NucDimSmale} that these $\mathrm{C}^*$-algebras are always classifiable using advanced techniques from the classification programme, rather than tools from dynamics. Thus it is once again difficult to glean information about the underlying Smale spaces. In addition, \cite{DS:NucDimSmale} showed that both the stable and unstable $\mathrm{C}^*$-algebra have finite nuclear dimension and are $\mathcal{Z}$-stable. However since these $\mathrm{C}^*$-algebras are nonunital, classification by the Elliott invariant is still out of reach, see \cite[Theorem 4.8 and Question 4.9]{DS:NucDimSmale}.

Here we turn our attention to group actions on Smale spaces and the corresponding actions induced on the associated $\mathrm{C}^*$-algebras. Such group actions have been studied at the Smale space level in various places. A number of explicit examples are discussed in Section \ref{exAutoSmaleSpace}. The specific case of shifts of finite type is already well-studied; see \cite{AdlKitMar:FinGrpOnSFT, BoyLinRud:AutoOFSFT, Boy:OpenProblems, Nas:TextileSys} among many others. Boyle's survey article \cite{Boy:OpenProblems} and the references within are a good place to get an introduction to this vast field, in particular see open problems 12 and 13 in \cite{Boy:OpenProblems}.

The goal of the present paper is to study how an action of a group on a Smale space relates to the induced action on its $\mathrm{C}^*$-algebras as well as to the associated crossed product.  Sufficient conditions for a group action on a general $\mathrm{C}^*$-algebra which ensure structural properties of the crossed product (simplicity, for example, or finite nuclear dimension, or $\mathcal{Z}$-stability) have been studied. These might be thought of as the noncommutative interpretation of a free action.  Of particular interest for this paper are when an action has finite Rokhlin dimension or is strongly outer. We establish sufficient conditions for when a group acting on a mixing Smale space induces an action on each of its associated $\mathrm{C}^*$-algebras with one of these properties. From there, we may deduce results about the crossed products.  

The Rokhlin property takes its motivation from the Rokhlin Lemma of ergodic theory. The Rokhlin Lemma says that a measure-preserving, aperiodic integer action can be approximated by cyclic shifts.  Connes successfully adapted this to the noncommutative setting of von Neumann algebras to classify automorphisms of the hyperfinite  II$_1$ factor up to outer conjugacy \cite{ MR0394228,  MR0448101}. Many generalisations followed in the von Neumann setting. Taking motivation from von Neumann algebras, the corresponding $\mathrm{C}^*$-algebraic theory gradually emerged, first for the restricted case of UHF algebras (the most straightforward $\mathrm{C}^*$-algebraic interpretation of a II$_1$ factor) \cite{MR647069, MR763780}, and more generally in the work of Kishimoto (see for example \cite{MR1344136, Kis:RohlinUHF, MR1660386}) and Izumi \cite{Izu:FreeRok1, Izu:FreeRok2}. The presence of the Rokhlin property for an action of a group on a $\mathrm{C}^*$-algebra allows that certain properties from the $\mathrm{C}^*$-algbera to the crossed product. However, since arbitrary $\mathrm{C}^*$-algebras need not have many projections, the Rokhlin property itself is often too strict. For this reason, weaker properties, via higher dimensional versions or tracial versions, emerged, see for examples \cite{HirWinZac:RokDim, Sza:RokDimZd, HSWW:RDFlow, Gar:RD, MR2441593, MR3347172, MR2854735}.  These are much more general but still allow for good preservation of properties in the crossed product. For further details on the Rokhlin property, Rokhlin dimension, and related $\mathrm{C}^*$-algebraic properties; see for example \cite[Introduction]{HirWinZac:RokDim}.

To the authors' knowledge, this is the first time general actions on Smale space $\mathrm{C}^*$-algebras have been studied. However, a number of special cases have been considered. In particular, the most obvious action on the stable, unstable, and homoclinic algebras is the one induced from $\varphi$. In \cite{PutSpi:Smale}, Putnam and Spielberg showed that if $(X, \varphi)$ is mixing, then $\salg \rtimes_{\varphi} \Z$ and $\ualg \rtimes_{\varphi} \Z$ are purely infinite. On the other hand, $\halg \rtimes_{\varphi} \Z$ is stably finite and in the special case of a shift of finite type, Holton \cite{Hol:RohSFT} showed that the shift map (that is, $\varphi$ in the case of a shift of finite type) leads to an action on $\halg$ with the Rokhlin property. Automorphisms of subshifts and their associated $\mathrm{C}^*$-algebras have been studied by a number of authors, see for example \cite{Mat:AutoShift} and references therein. Finally, in \cite{Sta:FinActTilAlg}, Starling studies certain finite group actions on tiling spaces and a $\mathrm{C}^*$-algebra associated to them which is related to the $\mathrm{C}^*$-algebras in the present paper through work of Anderson and Putnam \cite{AndPut:TopInvSubTilCstar}. This is discussed briefly in Example \ref{tilingEx}.

It is worth noting that an action at the Smale space level also induces an action on each of $\salg \rtimes_{\varphi} \Z$, $\ualg \rtimes_{\varphi} \Z$ and $\halg \rtimes_{\varphi} \Z$. Properties of the resulting crossed product algebras can be inferred from results in this paper by considering the action of the group generated by the original group and $\varphi$ on $\salg$, $\ualg$ and $\halg$. For example, if $\alpha$ is an automorphism of the given Smale space, then often (if $\alpha$ has finite order or is a power of $\varphi$ some care would be required) the study of induced actions on $\salg \rtimes_{\varphi} \Z$, $\ualg \rtimes_{\varphi} \Z$ and $\halg \rtimes_{\varphi} \Z$ amounts to studying the associated $\Z^2$-actions on $\salg$, $\ualg$ and $\halg$. 

The paper is structured as follows. In Section 1 we introduce Smale spaces. Section 2 provides examples of group actions on Smale spaces as well as a couple of short proofs about free actions and effective actions. We move on to $\mathrm{C}^*$-algebras in Section 3. 

In Section 4 we consider the induced action on the homoclinic $\mathrm{C}^*$-algebra, $\halg$. Two of our main results are included in this section: that actions of finite groups on mixing $(X, \varphi)$ which act freely on $X$ induce actions on $\halg$ with finite Rokhlin dimension and actions of discrete groups acting effectively on $X$ induce strongly outer actions on $\halg$. We also consider $\mathbb{Z}$-actions in the case that $(X, \varphi)$ is irreducible, but not necessarily mixing. In Section 5 we look at the induced actions on the stable and unstable $\mathrm{C}^*$-algebras of a Smale space, $\salg$ and $\ualg$. We show that we can use what we know about the induced action on $\halg$ to determine the behaviour on $\salg$ and $\ualg$ and consider the resulting properties of the crossed products.  Finally, in Section~6 we collect some results about $\mathcal{Z}$-stability, nuclear dimension, and classification of the crossed products.

In summary, the main results of the present paper are the following:

\begin{ntheorem}
Let $G$ be a discrete group acting effectively on a mixing Smale space $(X, \varphi)$. Then the induced action on $\halg$ is strongly outer.
\end{ntheorem}

When $G$ is a countable amenable group, then the fact that $\halg \rtimes_{\beta} G$ is $\mathcal{Z}$-stable follows from \cite[Theorem 1.1]{Sat:Amen}. However, since $\salg$ and $\ualg$ are nonunital, we cannot apply that result directly. Nevertheless, we show they are indeed $\mathcal{Z}$-stable.

\begin{ntheorem}
Suppose $G$ is a countable amenable group acting on a mixing Smale space, $(X, \varphi)$. Then, for the induced actions $\beta^{(S)} : G \to \Aut (\mathrm{C}^*(S))$ and $\beta^{(U)} : G \to \Aut (\mathrm{C}^*(U))$,  the crossed products $\mathrm{C}^*(S) \rtimes_{\beta^{(U)}} G$ and $\mathrm{C}^*(U) \rtimes_{\beta^{(U)}} G$ are both $\mathcal{Z}$-stable.
\end{ntheorem}

\begin{ntheorem}
Let $G$ be a residually finite group acting on an irreducible Smale space $(X, \varphi)$. Suppose the induced action of G on $\halg$ has Rokhlin dimension at most $d$. Then the induced actions of $G$ on $\mathrm{C}^*(S)$ and $\mathrm{C}^*(U)$ both have Rokhlin dimension at most $d$. The same statement holds for Rokhlin dimension with commuting towers.
\end{ntheorem}

%
%
%
%

\section{Preliminaries}

Smale spaces were defined by Ruelle~\cite{Rue:ThermForm}.  

\bn
\begin{ndefn}\!\cite[Section 7.1]{Rue:ThermForm} \label{DefSmaleSpace}
Let $(X, d)$ be a compact metric space and let $\varphi : X \to X$ be a homeomorphism. The dynamical system $(X, \varphi)$ is called a Smale space if there are two constants $\epsilon_X > 0$ and $0<\lambda_X<1$ and a map, called the bracket map,  
\[ [ \cdot, \cdot] : X \times X \to X \]
which is defined for $x, y \in X$ such that $d(x,y)< \epsilon_X$. The bracket map is required to satisify the following axioms: 
\begin{itemize}
\item[B1.] $\left[ x, x \right] = x$, 
\item[B2.] $\left[ x, [ y, z] \right] = [ x, z]$,
\item[B3.] $\left[ [ x, y], z \right] = [ x,z ]$,
\item[B4.] $\varphi[x, y] = [ \varphi(x), \varphi(y)]$;
\end{itemize}
for $x, y, z \in X$ whenever both sides in the above equations are defined. The system also satisfies
\begin{itemize}
\item[C1.] For $x,y \in X$ such that $[x,y]=y$, we have $d(\varphi(x),\varphi(y)) \leq \lambda_X d(x,y)$ and
\item[C2.] For $x,y \in X$ such that $[x,y]=x$, we have $d(\varphi^{-1}(x),\varphi^{-1}(y)) \leq \lambda_X d(x,y)$.
\end{itemize}
\end{ndefn}
\en

If the bracket map exists, it is unique (see for example \cite[page 7]{Put:C*Smale}). Also, in order to avoid certain trivial cases, we will always assume that $X$ is infinite.

\bn
\begin{ndefn} \label{1.2}
Suppose $(X, \varphi)$ is a Smale space, $x \in X$, and $0<\epsilon\le \epsilon_X$ and $Y, Z \subset X$ a subset of points. Then we define the following sets
\begin{enumerate}
\item $X^S(x, \varepsilon)  :=  \left\{ y \in X \mid d(x,y) < \varepsilon, [y,x]=x \right\},$
\item $X^U(x, \varepsilon)  :=  \left\{ y \in X \mid d(x,y) < \varepsilon, [x,y]=x \right\}, $
\item $X^S(x)  :=  \left\{ y\in X \mid \lim_{n \rightarrow + \infty} d(\varphi^n(x), \varphi^n(y)) =0 \right\}, $
\item $X^U(x)   :=  \left\{ y\in X \mid \lim_{n \rightarrow - \infty} d(\varphi^n(x), \varphi^n(y)) =0 \right\},$
\item $X^S(Z) := \cup_{x\in Z} X^S(x)$ and $X^U(Z):=\cup_{x\in Z} X^U(x)$, and
\item $X^H(Y, Z) := X^S(Y) \cap X^U(Z).$
\end{enumerate}
We say $x$ and $y$ are stably equivalent and write $x \sim_s y$ if $y \in X^S(x)$. Similarly, we say $x$ and $y$ are unstably equivalent, written $x \sim_u y$, if $y \in X^U(x)$. Points $x,y \in X$ are homoclinic if $y \in X^S(x) \cap X^U(x)$, meaning both $\lim_{n \rightarrow \infty} d(\varphi^n(x), \varphi^n(y)) =0$ and $\lim_{n \rightarrow -\infty} d(\varphi^n(x), \varphi^n(y)) =0$.
\end{ndefn}
\en

\bn \label{MixIrr}
A topological dynamical system $(X, \varphi)$ (not necessarily a Smale space) is called \emph{mixing} if, for every ordered pair of nonempty open sets $U, V \subset X$, there exists $N \in \mathbb{Z}_{>0}$ such that $\varphi^n(U) \cap V$ is nonempty for every $n \geq N$.

A Smale space $(X, \varphi)$ is \emph{irreducible} if for every ordered pair of nonempty open sets $U, V \subset X$ there exists $N \in \mathbb{Z}_{>0}$ such that $\varphi^N(U) \cap V$ is nonempty. In this case, Smale's decomposition theorem states that $X$ can be written as a disjoint union of finitely many clopen subspaces $X_1, \dots, X_N$ which are cyclically permuted by $\varphi$ and each system $(X_i, \varphi^N|_{X_i})$ is mixing. 
\en

\bn For $k \in \mathbb{N} \setminus \{0\}$, let 
\[ {\rm Per}_k(X, \varphi) = \{ x \in X \mid \varphi^k(x) = x \} \]
denote the set of points of period $k$, and let 
\[ {\rm Per}(X, \varphi)  = \cup_{k \in \mathbb{N} \setminus \{0\}} {\rm Per}_k(X, \varphi) \]
denote the set of all periodic points.

If $(X, \varphi)$ is irreducible, then the set ${\rm Per}(X, \varphi)$ is dense in $X$ and for each $k$, the set ${\rm Per}_k(X, \varphi)$ is finite  \cite{Rue:ThermForm}.
\en

\bn
\begin{ntheorem}\cite[Theorem 1]{RueSul:Currents} \label{BowenMeasure}
Given a mixing Smale space $(X, \varphi)$ there exists a unique $\varphi$-invariant, entropy-maximizing probability measure $\mu_X$ on $X$ which, for every $x \in X$ and $0< \epsilon < \epsilon_X$, can be written locally as a product measure supported on $X^{u}(x, \epsilon) \times X^{s}(x, \epsilon)$.
\end{ntheorem}
\en

The measure described in Theorem~\ref{BowenMeasure} is called the \emph{Bowen measure}, which in the case of a shift of finite type is the Parry measure, see for example \cite[Section 9.4]{LindMarcus:SDandCod}.

\section{Group actions on Smale spaces}
Let $G$ be a topological group. By an action of $G$ on a Smale space $(X, \varphi)$ we mean a continuous group homomorphism $G \to \homeo(X)$  such that 
\[ g \varphi(x) = \varphi (gx) \text{ for every } x \in X \text{ and } g \in G.\]
For a given Smale space $(X, \varphi)$, its automorphism group is defined to be
\[
{\rm Aut}(X, \varphi) := \{ \beta : X \rightarrow X \: | \: \beta \hbox{ is a homeomorphism and }\beta \circ \varphi = \varphi \circ \beta\}.
\]

Let $X$ be a compact metric space and $G$ a topological group. An action of group $G \to \homeo(X)$ is \emph{free} if, for every $x \in X$, we have that $g(x) = x$ if and only if $g = \id$. The action of $G$ is \emph{effective} (or \emph{faithful}), if for every $g \in G\setminus\{e\}$, there is an $x \in X$ such that $gx \neq x$. If $(X, \varphi)$ is a Smale space then we say the action is free, respectively effective, if the action on $X$ is free, respectively effective.

\subsection{Examples of automorphisms on Smale spaces} \label{exAutoSmaleSpace}
Group actions on Smale spaces are quite ubiquitous. Here we discuss four familiar classes of Smale spaces, three of which are treated in Putnam's work on Smale spaces and $\mathrm{C}^*$-algebras \cite{Put:C*Smale}, and give examples of their automorphisms. In each case, these actions will be free or effective and thus are covered by the results in Sections~\ref{haction} and \ref{IndActSandUSec}.

\bn
\begin{nexample}[Smale space automorphism] \label{alpha} Let $(X, \varphi)$ be a Smale space. Our first and most obvious example is the action induced by the homeomorphism $\varphi$ of $X$. Clearly $\varphi$ commutes with itself and hence defines an action of $\mathbb{Z}$ on $(X, \varphi)$. The action induced on the associated $\mathrm{C}^*$-algebras will prove to be a useful tool in the sequel.  Of course, powers of $\varphi$ are also automorphisms of $(X, \varphi)$.
\end{nexample}
\en

\bn
\begin{nexample}[Shifts of finite type]  \label{SFTexampleActions}
In the irreducible case, shifts of finite type are exactly the zero dimensional Smale spaces, though of course they were well known before Ruelle's work. The automorphism group and group actions on shifts of finite type have been studied by many authors, see for example \cite{AdlKitMar:FinGrpOnSFT, BoyLinRud:AutoOFSFT, Boy:OpenProblems, Nas:TextileSys} along with references therein. 

The full two-shift is the Smale space $(X, \varphi)$ where $X = \Sigma_{[2]}=\{ 0, 1\}^{\Z}$ with $\varphi= \sigma$ defined to be the left shift.  The automorphism group of the full two-shift is large: it contains, for example, every finite group and the free group on two generators \cite{BoyLinRud:AutoOFSFT}.  Here we highlight two automorphisms of finite order. The latter comes from \cite[Exercise 1.5.8]{LindMarcus:SDandCod}.
 
Define $\beta_1 : \Sigma_{[2]} \rightarrow \Sigma_{[2]}$ via $(a_n)_{n\in \Z} \mapsto (a_n+ 1$ mod$(2))_{n \in \Z}$ and $\beta_2 : \Sigma_{[2]} \rightarrow \Sigma_{[2]}$ via $(a_n)_{n\in \Z} \mapsto (b_n)_{n\in \Z}$ where 
\[ 
b_n := a_n + a_{n-1}( a_{n+1}+1)a_{n+2} \hbox{ mod}(2) .
\]
One can compute that $\beta_1^2$ and $\beta_2^2$ are both the identity. Hence they are order two automorphisms of the full two shift. The $\frac{\Z}{2\Z}$-action induced from $\beta_1$ is free while that of $\beta_2$ is effective but not free.
\end{nexample}
\en

\bn
\begin{nexample}[Solenoids] \label{solenoid} A solenoid is a Smale space obtained from a stationary inverse limit of a metric space equipped with a surjective continuous map that is subject to certain conditions such as those in \cite{Wil:ExpAtt} or in \cite{Wie:SmaleLimits}. 

Some prototypical examples are obtained as follows. Let $S^1 \subseteq \C$ (as the unit circle), $k\ge 2$ be an integer, and define $g: S^1 \rightarrow S^1$  by $g(z)=z^k$. Then $(X, \varphi)$ is the Smale space given by $X= \underleftarrow{\lim}(S^1, g)$ and the map $\varphi : X \rightarrow X$ defined as
\[
(z_n)_{n \in \N} \mapsto (g(z_0), g(z_1), g(z_2), \ldots ) = (g(z_0) , z_0, z_1, \ldots ).
\]

Let $\beta_1: X \rightarrow X$ be the map $(z_n)_{n\in \N} \mapsto (\bar{z}_n)_{n\in \N}$, where $\bar{z}$ denotes complex conjugate of $z$. For any $k$, $\beta_1$ defines an order two automorphism. There are other automorphisms. For example, if $g(z) = z^6$, then we have automorphisms defined by
\[
\beta_2( (z_n)_{n\in \N}) = (z^2_n)_{n\in \N} \quad \hbox{ and } \quad \beta_3( (z_n)_{n\in \N}) = (z^3_n)_{n\in \N}.
\]
The ${\Z}/{2\Z}$-action induced from $\beta_1$ is effective but not free. While the $\Z$-actions induced by $\beta_2$ and $\beta_3$ are each effective but not free. The automorphism group of these examples are known (see \cite[Proposition 1.6]{Tez:AutExpAtt} for details). Further results concerning the automorphism group of similar examples can also be found in \cite{Tez:AutExpAtt}. 
\end{nexample}
\en

\bn
\begin{nexample}[Hyperbolic toral automorphisms] \label{ToralAuto} Let $d\ge 2$ be an integer. A hyperbolic toral automorphism is a Smale space $(X, \varphi)$ were $X=\R^d/  \Z^d$ and $\varphi$ is induced by a $d \times d$ integer matrix $A$ with the following properties:
\begin{enumerate}
\item $|\det(A)| =1$;
\item no eigenvalue of $A$ has modulus one.
\end{enumerate}
A specific example when $d=2$ is
\[
 A = \left( \begin{array}{cc} 2 & 1 \\ 1 & 1 \end{array} \right).
\]
To obtain an automorphism of $(X, \varphi)$ one can take $B\in M_d(\Z)$ with det$(B)=\pm 1$ such that $AB=BA$. For example, again when $d=2$,
\[
B=\left( \begin{array}{cc} -1 & 0 \\ 0 & -1 \end{array} \right)
\]
defines an order two automorphism of any hyperbolic toral automorphism; the induced $\frac{\Z}{2\Z}$-action is effective. 

For larger $d$ there are more interesting automorphisms. An explicit example taken from \cite{Pol:GroPerComTorAuto} is
\[
 A = \left( \begin{array}{ccc} 1 & -1 & 0 \\ -1 & 2 & -1 \\ 0 & -1 & 2 \end{array} \right) \hbox{ and } B= \left( \begin{array}{ccc} 2 & 0 & -1 \\ 0 & 1 & 1 \\ -1 & 1 & 2 \end{array} \right)
\]
In this example, both $A$ and $B$ are hyperbolic. More on automorphisms of hyperbolic toral automorphisms can be found in \cite{BaaRob:SymTorAut} and \cite{Pol:GroPerComTorAuto} along with the references therein.
\end{nexample}
\en

\bn
\begin{nexample} [Substitution tilings] \label{tilingEx}
Many people have studied the dynamics of substitution tiling systems. In \cite[Chapter 2]{Con:noncomm}, Connes associates an AF algebra to the Penrose tiling by kites and darts; other $\mathrm{C}^*$-algebraic constructions were considered by Kellendonk  (see for example \cite{Kell:NCGTile}). It was Putnam and Anderson who showed that, under the assumptions that the substitution map $\omega$ is one-to-one, primitive and the tiling space $\Omega$ is of finite type (that is, has finite local complexity), that $(\Omega, \omega)$ is in fact a Smale space \cite{AndPut:TopInvSubTilCstar}. Moreover, the unstable $\mathrm{C}^*$-algebra associated to $(\Omega, \omega)$ is Morita equivalent to the tiling $\mathrm{C}^*$-algebra studied by Kellendonk \cite{AndPut:TopInvSubTilCstar}.

Many aperiodic tilings have interesting rotational symmetries. For example, the dihedral group, $D_5$ acts on the Penrose tiling \cite{SchBleGraCah:MetPhaOri}. In \cite{Sta:FinActTilAlg}, Starling studied free actions of finite subgroups of the symmetry group of these substitution tilings. The results of the present paper, in the special case of tilings, can be related to those in \cite{Sta:FinActTilAlg} by showing that the Morita equivalence considered by Anderson and Putnam can be made equivariant with respect to the given group action.
\end{nexample}
\en

Returning to general actions on Smale spaces, we begin by showing free actions are quite rare. However, as we shall see in Section~\ref{stronglyOuter}, an effective action is sufficiently strong to guarantee good properties of the induced action on the associated $\mathrm{C}^*$-algebras. The next two propositions are likely not new, but we could not find a reference for them, except for the case of a shift of finite type \cite{AdlKitMar:FinGrpOnSFT}. 

\bn
\begin{nprop} \label{notManyFreeActions}
Suppose $G$ has an element of infinite order and acts on a Smale space $(X, \varphi)$. Then $G$ does not act freely.
\end{nprop}

\begin{nproof}
Let $g\in G$ be an element of infinite order. The set periodic points of $(X, \varphi)$ is non-empty because the non-wandering set is non-empty by \cite[Appendix A.2]{Rue:ThermForm} and the periodic points are dense in the non-wandering set by \cite[Section 7.3]{Rue:ThermForm}. Let $x\in X$ be a periodic point with minimal period $k$. Then, for any $n \in \Z$,  
$$ \varphi^k (g^n  x) = g^n  \varphi^k (x) = g^n  x$$
Hence, for any $n \in Z$, $g^n x \in {\rm Per}_k(X, \varphi)$. Since ${\rm Per}_k (X, \varphi)$ is finite, there are $n_1 \neq n_2 \in \Z$ such that 
$$g^{n_1}  x = g^{n_2}x$$
It follows that the action is not free.
\end{nproof}
\en

\bn
\begin{nprop} \label{notEqualEffSet}
Suppose $G$ acts effectively on a mixing Smale space $(X, \varphi)$. For each $g\in G\setminus\{e\}$, the set
\[ 
\{ x \in X \: | \: gx \neq x \}
\]
is dense and open. Moreover, the Bowen measure of this set is one.
\end{nprop}

\begin{nproof}
The set $\{ x \in X \: | \: gx =x \}$ is closed since the map $G \times X \to X$ is continuous. Thus $\{ x \in X \: | \: gx \neq x \}$ is open.

Since $(X, \varphi)$ is mixing, there is a point $x_0 \in X$ with dense orbit. We show $x_0 \in \{ x \in X \: | \: gx \neq x \}$. Suppose not. Then for each $n\in \Z$,
\[ 
g(\varphi^n(x_0))= \varphi^n( g x_0) = \varphi^n(x_0).
\]
Thus, because the orbit of $x_0$ is dense, $g$ is the identity on $X$. This contradicts the assumption that $G$ acts effectively, so $x_0 \in \{ x \in X \: | \: gx \neq x \}$. Moreover, for each $n\in \Z$, $\varphi^n(x_0) \in  \{ x \in X \: | \: gx \neq x \}$. Since the orbit of $x_0$ is dense, $\{ x \in X \: | \: gx \neq x \}$ is also dense.

Finally, the set of points with dense orbit has Bowen measure one, which follows from Bowen's theorem and the fact that this is true for a shift of finite type \cite[Theorem 9.4.9]{LindMarcus:SDandCod}. The set $\{ x \in X \: | \: gx \neq x \}$ contains the set of points with dense orbit and hence also has Bowen measure one.
\end{nproof}
\en

\section{From Smale spaces to $\mathrm{C}^*$-algebras} \label{SmaSpaToCStar}

Ruelle was the first person to associate operator algebras to Smale spaces in \cite{Rue:NCAlgs}. We follow the approach introduced by Putnam and Spielberg \cite{Put:C*Smale, PutSpi:Smale}: three $\mathrm{C}^*$-algebras are constructed via the groupoid $\mathrm{C}^*$-algebra construction for  \'etale equivalence relations which capture the contracting, expanding, and asymptotic behaviour of the system given in Definition~\ref{1.2}.

%
%
Fix an irreducible Smale space $(X, \varphi)$. To define topologies on each of our equivalence relations, we first note the following. Using the notation from Definition~\ref{1.2}, it is not difficult to show that for any $0 < \epsilon \leq \epsilon_X$ we have
\[ X^U(x) = \cup_{n \in \mathbb{N}} \varphi^{n} (X^U(\varphi^{-n}(x), \epsilon)),\] and similarly that 
\[ X^S(x) = \cup_{n \in \mathbb{N}} \varphi^{-n}(X^S(\varphi^n(x), \epsilon)).\]
Each $\varphi^{n} (X^U(\varphi^{-n}(x), \epsilon))$, $n \in \mathbb{N}$ is given the relative topology from $X$ while $X^U(x)$ and $X^S(x)$ are given the topology coming from these inductive unions, see \cite[Page 10]{Put:C*Smale} for the precise details.
%
%

We could proceed to construct $\mathrm{C}^*$-algebras directly from the equivalence relations in Definition~\ref{1.2} following the construction of Putnam in \cite{Put:C*Smale}. However neither the stable nor the unstable groupoids would have a natural \'etale topology. Instead, following \cite{PutSpi:Smale}, we restrict our relation to those points equivalent to periodic points.

\bn
\begin{ndefn}
Let $P$ and $Q$ be finite $\varphi$-invariant sets of periodic points of $(X, \varphi)$. Define the stable and unstable groupoids of $(X, \varphi)$ by
\[ \G_S(P) := \{ (x,y) \in X^U(P) \times X^U(P) \mid x \sim_s y \},\]
and \[ \G_U(Q) := \{ (x,y) \in X^S(Q) \times X^S(Q) \mid x \sim_u y \}.\]
\end{ndefn}
\en
Up to Morita equivalence of groupoids these constructions do not depend on the choice of periodic points.

\bn
\begin{ndefn} Define the homoclinic groupoid of $(X, \varphi)$ by
\end{ndefn}
\[ \G_H := \{ (x, y) \in X \times X \mid x \sim_h y\}.\]
\en

Now, if $(v, w) \in X^S(P)$, then $v \sim_s w$ so there is some sufficiently large $N \in \mathbb{N}$ such that $d(\varphi^N(v), \varphi^N(w)) < \epsilon_X/2$. By continuity of $\varphi$, we may choose $\delta > 0$ small enough so that $\varphi^N(X^U(w , \delta)) \subset X^U(\varphi^N(w), \epsilon_X/2)$ and also $\varphi^N(X^U(v , \delta)) \subset X^U(\varphi^N(v), \epsilon_X/2)$. Then define
\[ h^s := h^s(v, w, N, \delta) : X^U(w, \delta) \to X^U(v, \epsilon_X), \quad x \mapsto \varphi^{-N}([\varphi^N(x), \varphi^N(v)]). \]
By \cite[Section 7.15]{Rue:ThermForm} this is a local homeomorphism. 

For any such $v, w, \delta, h , N$, we then define an open set by 
\[ V(v, w, \delta, h^s, N) := \{ (h^s(x), x) \mid x \in X^U(w, \delta) \} \subset \G_S(P).\] 
These sets generate an \'etale topology for $\G_S(P)$ \cite[Section 4]{PutSpi:Smale}. 

The construction for the topologies of $\G_U(Q)$ and $\G_H$ are similar; we refer the reader to \cite{PutSpi:Smale} for details. 

\subsection{$\mathrm{C}^*$-algebras}

Fix finite $\varphi$-invariant sets $P$ and $Q$ and let \[\mathcal{H} = \ell^2(X^H(P,Q))\]
where $X^H(P,Q)$ is the set of points in $X$ which are both stably equivalent to a point in $P$ and unstably equivalent to a point in $Q$. It is shown in \cite{Rue:ThermForm} that $X^H(P, Q)$ is countable.

If $\G$ is one of $\G_H, \G_S(P), \G_U(Q)$, let $C_c(\G)$ denote the a compactly supported functions on $\G$ with convolution product, 
\[ (f \ast g) (x, y) = \sum_{z \sim x, z \sim y} f(x,z) g(z, y), \quad (x,y) \in \G, \]
and 
\[f^*(x, y) = \overline{f(y,x)}, \quad (x,y) \in \G.\]

We can represent each of $C_c(\G_H)$, $C_c(\G_S(P))$ and $C_c(\G_U(Q))$ on the Hilbert space $\mathcal{H}$, and define the homoclinic algebra by
\[ \halg := \overline{C_c(\G_H)}^{\| \cdot \|_{\mathcal{H}}},\]
the stable algebra
\[ \salg := \overline{C_c(\G_S(P))}^{\| \cdot \|_{\mathcal{H}}},\]
and the unstable algebra
\[ \ualg := \overline{C_c(\G_U(Q))}^{\| \cdot \|_{\mathcal{H}}}.\]

\bn \begin{nremarks} \label{remRepSU}
1. We suppress the reference to $P$ and $Q$ in the notation of $\salg$ and $\ualg$. For any choice of such $P$ and $Q$, the resulting groupoids are Morita equivalent, and hence so are their $\mathrm{C}^*$-algebras. Thus from the perspective of most $\mathrm{C}^*$-algebraic properties, we don't need to keep track of the original choice.

2. In the usual groupoid $\mathrm{C}^*$-algebra construction for a groupoid $\mathcal{G}$, the algebra of compactly supported functions $C_c(\mathcal{G})$ is represented on the Hilbert space $\ell^2(\mathcal{G})$ and the completion is the reduced groupoid $\mathrm{C}^*$-algebra $\mathrm{C}_r(\mathcal{G})$. However, when the groupoid is amenable, the completion of any faithful representation will result in the same $\mathrm{C}^*$-algebra. Here $\G_H$, $\G_S(P)$ and $\G_U(Q)$ are each amenable \cite{PutSpi:Smale}. It is convenient, however, to represent them all on the same Hilbert space (namely $\ell^2(X^H(P,Q))$) because we can consider the interactions between operators coming from the different algebras. This will be particularly useful when showing how finite Rokhlin dimension passes from $\halg$ to $\salg$ and $\ualg$ in Section~\ref{IndActSandUSec}.
\end{nremarks}
\en

\bn \label{knownprops}
Let $(X, \varphi)$ be a mixing Smale space. Then $\halg$, $\salg$ and $\ualg$ are simple by \cite{Put:C*Smale} and \cite[Theorem 1.3]{PutSpi:Smale}, separable (since $X$ is a metric space) and nuclear by \cite[Theorem 1.3]{PutSpi:Smale} and \cite[Theorem 3.1]{Put:C*Smale}. The homoclinic $\mathrm{C}^*$-algebra $\halg$ is unital since the diagonal $X \times X$ is open in $\G_H$ and $X$ is compact.


Each of $\halg$, $\salg$ and $\ualg$ admit a trace \cite[Theorem 3.3]{Put:C*Smale}, hence are stably finite.  The traces on $\salg$ and $\ualg$ are not bounded while $\halg$ admits a tracial state. Moreover, when $(X, \varphi)$ is mixing, this trace is unique \cite{Hou:SmaleTrace}. 
  
 In  \cite{DS:NucDimSmale}, the authors showed that for a mixing Smale space, $\halg$, $\salg$ and $\ualg$ have finite nuclear dimension and hence are $\mathcal{Z}$-stable. (There it was not noted that $\salg$ and $\ualg$ are $\mathcal{Z}$-stable; however this follows from \cite[Corollary 8.7]{Tik:Non1Zstab}. In fact, it can be proved directly in a similar manner to Theorem~\ref{zstab} below, once we know that $\halg$ is $\mathcal{Z}$-stable.) It then follows that $\halg \otimes \uhf$ is tracially approximately finite (TAF) in the sense of \cite{Lin:TAF1} for any UHF algebra of infinite type $\uhf$. In particular, the \mbox{$\mathrm{C}^*$-algebras} coming from the homoclinic relation on a mixing Smale space is classified by the Elliott invariant, see \cite[Theorem 4.7]{DS:NucDimSmale}.
 \en

Suppose $G$ is a discrete group acting on a mixing Smale space. For the action to induce a well-defined action on $\G_U$ and $\G_S$, the choice of finite sets of $\varphi$-invariant periodic points must be $G$-invariant. Fortunately, this can always be arranged.

\bn
\begin{nlemma}
Let $G$ be a discrete group acting effectively on a mixing Smale space $(X , \varphi)$. Then there exists a finite set of $\varphi$-invariant periodic points $P$ such that $gp \in P$ for every $g \in G$ and every $p \in P$.
\end{nlemma}

\begin{nproof}
Let $P'$ be any finite $\varphi$-invariant set of periodic points. Then $P' \subseteq {\rm Per}_n(X, \varphi)$ for some $n\in \N$. We know that ${\rm Per}_n(X, \varphi)$ is finite. Also, for any $g\in G$ and $x\in {\rm Per}_n(X, \varphi)$ it was shown in the proof of Proposition~\ref{notManyFreeActions} that $gx \in {\rm Per}_n(X, \varphi)$. It follows that the set 
\[
P=\{ p \in X \mid p=gx \hbox{ for some }g\in G, x\in P'\} 
\]
is contained in ${\rm Per}_n(X, \varphi)$ and hence is finite. It is $G$-invariant by construction and it is $\varphi$-invariant because $g \varphi (x) = \varphi( gx)$ for any $g\in G$ and $x\in X$.
\end{nproof}
\en

For the remainder of the paper, we will assume that $P$ and $Q$ are $G$-invariant. 

Let  $\G$ be one of $\G_H, \G_S(P), \G_U(Q)$ and let $G$ be a group acting on $(X, \varphi)$. Since $P$ and $Q$ are assumed to be $G$-invariant,   $g(x,y) \mapsto (gx, gy)$  defines an induced action of $G$ on $\G$.  The action of $G$ on $\G$ in turn induces an action on $\mathrm{C}^*(\G)$.

\bn
\begin{nexample} \label{alphaCstar}
Following Example~\ref{alpha}, the Smale space homeomorphism $\varphi$ induces a $\mathbb{Z}$-action on each of $\halg, \salg , \ualg$. These actions are denoted by $\alpha$, $\alpha_S$ and $\alpha_U$, respectively. Properties of $\alpha$ (in particular, \cite[Theorem 3.2]{Put:C*Smale}) will prove indispensable for the results in the next section and also in Section~\ref{IndActSandUSec}, where we seek to pass from known properties about $\halg$ to the nonunital $\salg$ and $\ualg$.
\end{nexample}
\en

\bn
\begin{nexample}
Suppose $(\Sigma, \sigma)$ is a mixing shift of finite type. Then the $C^*$-algebras associated to $\G_H$, $\G_S(P)$ and $\G_U(Q)$ are each AF. Moreover if $G$ is a finite group acting effectively on $(\Sigma, \sigma)$, then using \cite[Observations 1 and 2]{AdlKitMar:FinGrpOnSFT} one can show that the induced action of $G$ on each of these AF $\mathrm{C}^*$-algebras is locally representable in the sense of \cite{HandRoss:CompAF}. This is the case for the automorphisms $\beta_1$ and $\beta_2$ on the full two-shift given in Example \ref{SFTexampleActions}.  It follows from results in \cite{HandRoss:CompAF} that the crossed products associated to the action of $G$ are each also AF. 
\end{nexample}
\en

\bn
\begin{nexample}
Suppose $(X, \varphi)$ is the Smale space obtained via the solenoid construction in Example \ref{solenoid} with $Y=S^1$ and $g(z)=z^n$ ($n\ge 2$). Then the stable and unstable algebras are the stabilisation of a Bunce--Deddens algebra. For details in the case $n=2$, see page 28 of \cite{Put:C*Smale}. Automorphisms of Bunce--Deddens algebras are considered in \cite{Pas:HomBDalg}, for example. 
\end{nexample}
\en

\bn
\begin{nexample}
If $(X, \varphi)$ is a hyperbolic toral automorphism, as in Example \ref{ToralAuto}, then the stable and unstable algebras are the stabilization of irrational rotation algebras as is shown on pages 27-28 of \cite{Put:C*Smale}. Automorphisms of irrational rotation algebras are well-studied, see \cite{EchLueckPhiWal:irrrat} and reference therein.
\end{nexample}
\en

\section{The induced action on $\halg$} \label{haction}

In the sequel, the aim is to provide conditions of a group action on a mixing Smale space which will allow us to determine structural properties of the crossed products of the associated $\mathrm{C}^*$-algebras by the induced group action. The idea is to determine what properties are preserved when passing from the $\mathrm{C}^*$-algebra to its crossed product. If we interpret a $\mathrm{C}^*$-crossed product as a ``noncommutative orbit space'' then what we are asking for is some sort of \,  ``freeness'' condition. In the $\mathrm{C}^*$-algebraic context this might take a number of different forms. Here, we focus on the Rokhlin dimension of an action (which is akin to a ``coloured'' version of noncommutative freeness), initially proposed for finite group and $\mathbb{Z}$-actions by Hirshberg, Winter and Zacharias \cite{HirWinZac:RokDim} and subsequently generalised to other groups \cite{Sza:RokDimZd, Gar:RD, SzaWuZac:RP, HSWW:RDFlow, BroTikZel:RD}, as well as the notion of a ``strongly outer action'' (which can be thought of as a noncommutative approximation of freeness in trace) (see for example, \cite{MatSat:StrongOuter}). Recent results of Sato \cite{Sat:Amen} also play a key role, although for finite group actions and $\Z^d$-actions these new results are not required. 

\subsection{Finite group actions on  $\halg$}

For this subsection, we fix a mixing Smale space $(X, \varphi)$. Let $\halg$ denotes its homoclinic algebra. Here we study actions induced on $\halg$ by free actions of $G$ on $(X, \varphi)$. Given an action of a group $G$ on a $C^*$-algebra, $A$, $\beta : G \rightarrow {\rm Aut}(A)$ we will write $\beta_h$ for $\beta(h)$. Based on Proposition \ref{notManyFreeActions} which shows that freeness is unlikely to hold for infinite groups, we restrict to the case that $G$ is finite. 

\bn
\begin{ndefn}\!\!\cite[Definition 1.1]{HirWinZac:RokDim} \label{RokDimFinGrp}
Let $A$ be a unital \mbox{$\mathrm{C}^*$-algebra} and let $G$ be a finite group. An action $\beta : G \to \Aut(A)$ has Rokhlin dimension $d$ if $d$ is the least integer such that the following holds: for each $\epsilon>0$ and each finite subset $F \subset A$ there are positive contractions
$$ \left( f^{(l)}_g \right)_{l=0, \ldots d; g\in G} \subset A $$
such that
\begin{enumerate}
\item $\| f^{(l)}_g f^{(l)}_h \| < \epsilon$, for all $l \in \{ 0, \ldots , d\}$ $g \neq h$ in $G$;
\item $\| \sum_{l=0}^d \sum_{g \in G} f^{(l)}_g - 1 \| < \epsilon$;
\item $\| \beta_h(f^{(l)}_g) - f^{(l)}_{hg} \| < \epsilon$  for all $l \in \{ 0, \ldots , d\}$, $g, h \in G$;
\item $\| [ f^{( l )}_g, a] \| < \epsilon$ for all $l \in \{ 0, \ldots , d\}$, $g \in G$ and $a \in F$.
\end{enumerate}
\end{ndefn}
\en

When $d = 0$ the action is said to have the Rokhlin property. In this case the contractions $(f_g)_{g \in G}$ can in fact be taken to be projections.  The definition of the Rokhlin property for $\mathrm{C}^*$-algebras was introduced by Izumi \cite{Izu:FreeRok1, Izu:FreeRok2}.

\bn
\begin{nlemma} \label{lastConIsClear}
Let $G$ be a finite group acting on $(X, \varphi)$ and denote by  $\beta : G \to \Aut(\halg)$ the induced action. Let $d$ be a nonnegative integer. Suppose that for any $\epsilon>0$ there are positive contractions
$$ \left( f^{(l)}_g \right)_{l=0, \ldots d; g\in G} \subset \halg $$
such that
\begin{enumerate}
\item $\| f^{(l)}_g f^{(l)}_h \| < \epsilon$, for $g \neq h \in G$, $l \in \{0, \dots, d\}$;
\item $\| \sum_{l=0}^d \sum_{g \in G} f^{(l)}_g - 1 \| < \epsilon$;
\item $\| \beta_h(f^{(l)}_g) - f^{(l)}_{hg} \| < \epsilon$  for all $l \in \{ 0, \ldots , d\}$, $g, h \in G$.
\end{enumerate}
Then $\beta$ has Rokhlin dimension at most $d$.
\end{nlemma}

\begin{nproof}
Let $\epsilon>0$ and $F \subset \halg$ be a finite set. Take positive contractions 
$$ \left( f^{(l)}_g \right)_{l=0, \ldots d; g\in G} \subset \halg $$
with the properties assumed in the statement of the theorem.

By \cite[Theorem 3.2]{Put:C*Smale}, there is $n \in \N$ such that $\| [ \alpha^n(  f^{( l )}_g) , a] \| < \epsilon$ for all $l \in \{ 0, \ldots , d\}$, $g \in G$ and $a \in F$, where $\alpha$ is the automorphism from Example~\ref{alphaCstar}.  Since $\alpha$ is an automorphism, 
$$\left( \alpha^n( f^{(l)}_g) \right)_{l=0, \ldots d; g\in G} \subset \halg $$
satisfies the other requirements of Definition \ref{RokDimFinGrp}.
\end{nproof}
\en

\bn
\begin{ncor} \label{FinFreActCor}
Suppose $A$ is a $G$-invariant unital $\mathrm{C}^*$-subalgebra of $\halg$. If $G$ acting on $A$ has Rokhlin dimension at most $d$ then the action of $G$ on $\halg$ also has Rokhlin dimension at most $d$. In particular, if $G$ acting on $C(X)$ has Rokhlin dimension at most $d$, then $G$ acting on $\halg$ also has Rokhlin dimension at most $d$.
\end{ncor}

\begin{nproof}
By assumption, given $\epsilon>0$, there exists  positive contractions 
$$ \left( f^{(l)}_g \right)_{l=0, \ldots d; g\in G} \subset A \subset \halg $$
such that the hypotheses of Lemma \ref{lastConIsClear} hold; this then implies the result. 

When considering the statement concerning $C(X)$ one need only note that $C(X)$ is a $\mathrm{C}^*$-subalgebra of $\halg$ and that it is invariant under the action of $G$.
\end{nproof}
\en

\bn
\begin{ncor}
Suppose $(X, \varphi)$ is a mixing Smale space and a finite group $G$ acts on $(X, \varphi)$ freely. Then $G$ acting on $\halg$ has finite Rokhlin dimension.
\end{ncor}

\begin{nproof}
For finite group actions on a compact space, freeness implies finite Rokhlin dimension \cite[Proposition 2.11]{HirPhi:RD}. The result then follows from Corollary \ref{FinFreActCor}.
\end{nproof}
\en

\bn
\begin{ncor}
Suppose $(\Sigma, \sigma)$ is a mixing shift of finite type and a finite group $G$ acts on $(\Sigma, \sigma)$ freely. Then the action of $G$  on $\halg$ has the Rokhlin property.
\end{ncor}

\begin{nproof}
Since $\Sigma$ is the Cantor set and $G$ acts freely, the action of $G$ on $C(\Sigma)$ has the Rokhlin property. Corollary \ref{FinFreActCor} then implies the result.
\end{nproof}
\en

\subsection{Strongly outer actions} \label{stronglyOuter}
Using \cite[Remark 2.8]{MatSat:StrongOuter}, we prove the first of the three results listed at the end of the introduction; it  appears as Theorem \ref{SO}. To begin, we recall the Vitali covering theorem and the definitions needed for its statement.

\bn 
\begin{ndefn}
A finite measure $\mu$ on a metric space $(X, d)$ is said to be \emph{doubling} if there exists a constant $M > 0$ such that 
\[ \mu(B(x, 2\epsilon)) \leq M \mu(B(x, \epsilon)) \]
for any $x \in X$ and any $\epsilon >0$.
\end{ndefn}
\en

\bn
\begin{ndefn}
Suppose $(Y,d)$ is a metric space and $A \subseteq Y$. Then a Vitali cover of $A$ is a collection of closed balls $\mathcal{B}$ such that inf$\{ r>0 \mid B(x,r) \in \mathcal{B} \} =0$ for all $x\in A$.
\end{ndefn}
\en

\bn
\begin{ntheorem}\textup{[Vitali Covering Theorem]}
Suppose $(Y,d)$ is a compact metric space, $\mu$ is a doubling measure, $A \subseteq Y$ and $\mathcal{F}$ is Vitali cover of $A$. Then, for any $\epsilon>0$, there exists finite disjoint family $\{ F_1, F_2, \ldots F_n \} \subseteq \mathcal{F}$ such that $\mu( A  - \cup^n_{i=1} F_i) < \epsilon$.
\end{ntheorem}
\en

\bn
\begin{nlemma} \label{weakRPish}
Suppose $(X, \varphi)$ is a mixing Smale space with Bowen measure $\mu$ and $G$ is a discrete group acting effectively on $(X, \varphi)$. Let $\beta : G \to \Aut(\halg)$ denote the induced action. Then, for any $g\in G \setminus \{e\}$ and $\epsilon>0$, there exists positive contractions $(f_i)_{i=1}^k \subset \halg$ such that
\begin{enumerate}
\item $f_i \cdot f_j=0$ for all $i\neq j$;
\item $f_i \cdot \beta_g(f_i) =0$;
\item $\tau(\sum_{i=1}^k f_i)> 1- \epsilon$ where $\tau$ denotes the (unique) trace on $\halg$ obtained from $\mu$.
\end{enumerate}
\end{nlemma}

\begin{nproof} 
Fix $g\in G \setminus \{e\}$ and $\epsilon>0$. Let $D=\{ x \in X \: | \: gx \neq x \}$. Then $D$ is open and has full measure by Proposition \ref{notEqualEffSet}. Moreover, for each $x\in D$ there exists $\delta_x>0$ such for any $0<\delta\le\delta_x$ we have $\overline{B_{\delta}(x)} \cap g(\overline{ B_{\delta}(x)}) = \emptyset$ and $\overline{B_{\delta}(x)} \subseteq D$.  

Let 
$$\mathcal{F}= \left\{ F_x \: | \: F_x = \overline{ B_{\delta}(x)} \hbox{ for some }x\in D \hbox{ and some }0<\delta\le \frac{\delta_x}{2} \right\}.$$
By construction, the collection $\mathcal{F}$ is a Vitali covering of $D$. The Bowen measure is doubling (see for example \cite[Introduction]{PesWei:FracAnaConExpMap}) so we may apply the Vitali Covering Theorem to obtain a finite subcollection of $\mathcal{F}$, $\{ F_{x_1}, \ldots, F_{x_k} \}$, with the following properties
\begin{enumerate}
\item $F_{x_i} \cap F_{x_j}=\emptyset$;
\item $\mu(\cup_{i=1}^k F_{x_i}) > 1- \epsilon$.
\end{enumerate}
Recall that $\mu(D)=1$.

Since each $F_{x_i}$ is compact (they are closed in a compact space) and pairwise disjoint, we can use Urysohn's lemma and the Tietze extension theorem to obtain pairwise disjoint open sets $(U_i)_{i=1}^k$ such that for each $i$, $F_{x_i}\subseteq U_{x_i}\subseteq B_{\delta_{x_i}}(x_i)$ and functions $(f_i)_{i=1}^k \subseteq C(X) \subseteq \halg$ with the following properties:
\begin{enumerate}
\item $0\le f_i \le 1$,
\item $ {\rm supp}(f_i) \subseteq U_i$,
\item $F_{x_i} \subseteq \{ x \: |\: f_i(x)=1 \}$;
\end{enumerate}
for every $i \in \{1, \dots, k\}$.

Finally, we show that $(f_i)_{i=1}^k$ has the required properties. They are by definition positive contractions. Moreover,
\begin{enumerate}
\item $f_i \cdot f_j =0$ for every $i \neq j \in \{1, \dots, k\}$, since ${\rm supp}(f_i) \subseteq U_i$ and $U_{i} \cap U_{j}=\emptyset$;
\item $f_i \cdot \beta_g(f_i)=0$ for every $i \in \{0, \dots, k\}$ since $U_i \subseteq B_{\delta_{x_i}}(x_i)$ implies that $U_i \cap g( U_i) = \emptyset$.
\end{enumerate}
Finally,
\begin{eqnarray*}
 \tau\left( \sum_{i=1}^k f_i \right) &  \ge & \sum_{i=1}^k \mu( \{ x | f_i(x) =1 \} ) \\
& \ge & \sum_{i=1}^k \mu(F_i)  \\
& = & \mu( \cup_{i=1}^k F_i )  \\
& > & 1 - \epsilon,
\end{eqnarray*}
showing (iii) holds.
\end{nproof}
\en

\bn Let $A$ be a $\mathrm{C}^*$-algebra and $\tau$ a state on $A$. We denote by $\pi_{\tau}$ the representation of $A$ corresponding to the GNS construction with respect to $\tau$. In this case, $\pi_{\tau}(A)''$ is the enveloping von Neumann algebra of $\pi_{\tau}(A)$.

\begin{ndefn}\!\!\cite[Definition 2.7 (2), (3)]{MatSat:StrongOuter} Let $A$ be a unital simple $\mathrm{C}^*$-algebra with nonempty tracial state space $T(A)$. An automorphism $\beta$ of $A$ is \emph{not weakly inner} if, for every $\tau \in  T(A)$ such that $\tau \circ \beta = \tau$,  the weak extension of $\beta$ to $\pi_{\tau}(A)''$ is outer. If $G$ is a discrete group, then an action $\beta: G \to \Aut(A)$ is \emph{strongly outer} if, for every $g \in G \setminus \{e\}$, the automorphism $\beta_g$ is not weakly inner.\end{ndefn}
\en

The next lemma and theorem are based on arguments due to Kishimoto in \cite{Kis:RohlinUHF}.

\bn
\begin{nlemma} \label{easy}
Let $A$ be a unital $\mathrm{C}^*$-algebra and $\tau \in T(A)$.  Then for every $\epsilon > 0$ there is $\delta > 0$ such that for any positive contraction $f \in A$ such that $\tau(f) > 1 - \delta$ we have $\tau(a) \leq \tau(fa)  + \epsilon$ for every $a \in A$.
\end{nlemma}

\begin{nproof} It is enough to show this holds when $a \in A_+$. Given $\epsilon > 0$ let $\delta = \epsilon^2$. Then since $a = af + a(1-f)$ we have
\begin{eqnarray*}
\tau(a) - \tau(af) &=& \tau(a(1-f)) \\
&\leq&  \tau(a^2)^{1/2}\tau((1-f)^2)^{1/2} \\
&\leq& \tau((1-f)^2)^{1/2}\\
&\leq& \tau(1-f)^{1/2} \\
&\leq& \epsilon. 
\end{eqnarray*} 
\end{nproof}
\en

\bn
\begin{ntheorem} \label{SO} Let $(X, \varphi)$ be a mixing Smale space with an effective action of a discrete group $G$. Then the induced action $\beta : G \to \Aut(\halg)$ is strongly outer.
\end{ntheorem}

\begin{nproof}
Let $\tau$ denote the tracial state on $\halg$ corresponding to the Bowen measure on $(X, \varphi)$. Fix $g \in G\setminus \{e\}$. Since $\tau$ is the unique tracial state, we have $\tau = \tau \circ \beta_g$. Let $\epsilon > 0$ and let $F \subset A$ be a finite subset. Let $\delta$ be the $\delta$ of Lemma~\ref{easy} with respect to $\epsilon/2$.   By Lemma~\ref{weakRPish}, we can find positive contractions $f_1, \dots, f_k \in \halg$ such that 
\begin{enumerate}
\item $f_if_j = 0$ for every $0 \leq i \neq j \leq k$,
\item $f_i \beta_g(f_i) = 0$, 
\item $\tau(f_1 + \dots + f_k) > 1 - \delta$. 
\end{enumerate}
There is a sufficiently large $N \in \mathbb{N}$ such that, for any $i =1, \dots, k$, we have  
\[   \| \alpha^N(f_i) a - a \alpha^N(f_i) \| < \epsilon/2 \text{, for any } a \in F.\]
Since $\alpha^N : \halg \to \halg$ is an automorphism, we have that
\begin{enumerate}[resume]
\item $\alpha^N(f_i) \alpha^N(f_j) = 0$ for every $1 \leq i \neq j \leq k$ and
\item $\tau(\alpha^N(f_1) + \dots + \alpha^N(f_k)) > 1 - \delta.$
\end{enumerate}
Moreover, since the action of $G$ on $X$ commutes with $\varphi$, we have, for any $f\in \halg$, that $\alpha^N(\beta_g(f)) = \beta_g(\alpha^N(f))$. Hence
\begin{enumerate}[resume]
\item $\alpha^N(f_i) \beta_g(\alpha^N(f_i)) = 0.$ \label{t2}
\end{enumerate}
Thus we have  
\[ \| \alpha^N(f_i) a \beta_g(\alpha^N(f_i)) \| < \epsilon/2\]
for any $a \in F$.

Consider the crossed product $\halg \rtimes_{\beta_g} \mathbb{Z}$. It is generated by $\halg$ and $u_g$ where $u_g$ is unitary satisfying $u_g a u_g^* = \beta_g(a)$ for every $a \in \halg$.  Let $\sigma \in T(\halg \rtimes_{\beta_g} \mathbb{Z})$ and note that $\sigma|_{\halg} = \tau$.  Put $f = \sum_{i=1}^k \alpha^N(f_i)$. By the above and Lemma~\ref{easy} we have
\begin{eqnarray*}
\sigma(a u_g) &\leq& \sigma(f a u_g f) +  \epsilon\\
&=&  \sum_{i=1}^k \sigma( \alpha^N(f_i) a \beta_g(\alpha^N(f_i)) u_g ) + \epsilon \\
&<& 2\epsilon.
\end{eqnarray*}
Since $\epsilon$ and $F$ were arbitrary, it follows that $\sigma(a u_g) = 0$ for all $a \in A$.  That $\beta_g$ is not weakly inner now follows from the proof of \cite[Lemma 4.4]{Kis:RohlinUHF}. Since this holds for every $g \in G\setminus \{e\}$, it follows that the action $\beta : G \to \Aut(\halg)$ is strongly outer.
\end{nproof}
\en

\bn
\begin{nprop} \label{ZdActRok} Suppose $(X, \varphi)$ is a mixing Smale space with an effective action of $\Z^m$. Then the induced action on $\halg$ has finite Rokhlin dimension. In fact, the Rokhlin dimension of the action is bounded by $4^m -1$. If $m = 1$ then the action has Rokhlin dimension no more than one.
\end{nprop}

\begin{nproof}
This result follows directly from \cite[Theorem 1.1]{Liao:RDimZm}.
\end{nproof}
\en

\subsection{$\mathbb{Z}$-actions on irreducible Smale spaces}

Throughout this subsection, we fix an irreducible Smale space $(X, \varphi)$ with decomposition into mixing components given by $X = X_1 \sqcup \dots \sqcup X_N$, as in \ref{MixIrr}.

To begin, we recall the definition of Rokhlin dimension for integer actions as given in \cite{HirWinZac:RokDim}.

\bn
\begin{ndefn}\cite[Definition 2.3]{HirWinZac:RokDim} \label{ZRD}
Let $A$ be a unital $\mathrm{C}^*$-algebra. An action of the integers $\beta : \mathbb{Z} \to \Aut(A)$ has Rokhlin dimension $d$ if $d$ is the least natural number such that the following holds: for any finite subset $\mathcal{F} \subset A$, and $p \in \mathbb{N}$ and any $\epsilon > 0$ there are positive contractions 
\[ f^{(l)}_{0,0}, \dots, f^{(l)}_{0, p-1}, f^{(l)}_{1, 0}, \dots, f^{(l)}_{1, p}, \qquad l \in \{0, \dots, d\},\]
in $A$ satisfying
\begin{enumerate}
\item for any $l \in \{0, \dots, d\}$ we have $\| f^{(l)}_{r,i} f^{(l)}_{s, j} \| < \epsilon$ whenever $(r,i) \neq (s, j)$,
\item $\| \sum_{l=0}^d \sum_{r =0}^1 \sum_{j = 0}^{p-1+r} f^{(l)}_{r,j} -1\| < \epsilon$,
\item $\|\beta_1(f^{(l)}_{r, j}) - f^{(l)}_{r, j+1} \| < \epsilon$ for every $r \in \{0,1\}$, $j \in \{0, \dots, p-2+r\}$ and $l \in \{0, \dots, d\}$, 
\item $\|\beta_1(f^{(l)}_{0,p-1} + f^{(l)}_{1, p}) - (f^{(l)}_{0,0} + f^{(l)}_{1,0}) \| < \epsilon$ for every $l \in \{0, \dots, d\}$
\item $\|[ f^{(l)}_{r,j}, a ] \| < \epsilon$ for every $r, j, l$ and $a \in \mathcal{F}$.
\end{enumerate}
\end{ndefn}
\en

As in the case of a finite group, in the special case of the homoclinic algebra, we don't need to worry about satisfying the last condition. The proof is similar to Lemma \ref{lastConIsClear} and is omitted.

\bn
\begin{nlemma} \label{lastConIsClearZ}
Let $\beta : \mathbb{Z} \to \Aut(\halg)$ be an action of the integers. Suppose that for any $p \in \mathbb{N}$ and $\epsilon > 0$ there are positive contractions 
\[ f^{(l)}_{0,0}, \dots, f^{(l)}_{0, p-1}, f^{(l)}_{1, 0}, \dots, f^{(l)}_{1, p}, \qquad l \in \{0, \dots, d\},\]
satisfying (i)--(iv) of Definition~\ref{ZRD}. Then $\beta$ has Rokhlin at most dimension $d$.
\end{nlemma}
\en

\bn
\begin{nlemma} \label{sigma} Let $\beta \in \Aut(X, \varphi)$. Then there is a $\sigma \in S_N$ such that 
\[ \beta|_{X_i} : X_i \to X_{\sigma(i)} \]
for each $i = 1,\dots, N$.
\end{nlemma}

\begin{nproof}
Since $(X_i, \varphi^N|_{X_i})$ is mixing, there is a point, call it $x_i$, with dense $\varphi^N$-orbit in $X_i$. Since $X_1, \dots, X_N$ are disjoint, there is $k(i) \in \{1, \dots, N\}$ such that $\beta(x_i) \in X_{k(i)}$ and $\beta(x_i) \notin X_k$ for $k \neq k(i)$. Now $\varphi^N \circ \beta(x_i) \in X_{k(i)}$, so  we have $\varphi^N \circ \beta(x_i) = \beta \circ \varphi^N(x_i) \in X_{k(i)}$. The fact that $x_i$ has dense $\varphi^N$-orbit then implies $\beta|_{X_i} : X_i \to X_{k(i)}$. 

Suppose that $\beta|_{X_j}(X_j) \subset X_{k(i)}$. If $j \neq i$, there exists some $l \in \{1, \dots, N\}$ such that $X_l$ is not in the image of any $\beta(X_i)$, $i \in \{ 1, \dots, N\}$. Choose $x \in X_l$. Then $\beta(x )\in X_{k(l)}$ and there is some $m$, $0<m<N$ satisfying $\varphi^m(X_{k(l)}) = X_{l}$.  But this implies $\beta \circ \varphi^m(x) = \varphi^m \circ \beta(x) \in X_l$. So we have $j= i$, which proves the theorem.
\end{nproof}
\en

\bn
\begin{nlemma} \label{LS}
Let $\beta \in \Aut(X, \varphi)$. Then there are  $L, S \in \mathbb{N}\setminus\{0\}$ such that $N = LS$ and partition of $\{1, \dots, N\}$ into $S$ subsets of $L$ elements $\{r_{s,1}, \dots, r_{s,L}\}$, $1 \leq s \leq S$ satisfying, for each $s$,
\[ \beta(X_{r_{s, l}}) = X_{s, r_{l+1 \text{mod} L}}.\]
\end{nlemma}

\begin{nproof}
Let $r_{1,1} = 1$. By the previous lemma there is $r_{1,2} \in \{1, \dots, N\}$ such that $\beta : X_{r_{1,1}} \to X_{r_{1,2}}$. If $r_{1, 2} \neq 1$ then there is $r_{1, 3} \neq r_{1,2}$ such that $\beta: X_{r_{1,2}} \to X_{r_{1,3}}$. By induction and the pigeonhole principal, there is some $0 < L \leq N$ such that $r_{1, l}$ are all distinct, $\beta: X_{r_{1,l}} \to X_{r_{1,l+1}}$, for $1 \leq l \leq L-1$ and $\beta  : X_{r_{1,L}} \to X_{r_{1,1}}$. If $L= N$, we are done.

Otherwise, there is some $X_{r_{2,1}}  \neq X_{r_{2, l}}$ for every $1 \leq l \leq L$. Arguing as above and using the previous lemma, there is some $L' \in \{1, \dots, N-L\}$ such that $r_{2, l}$ are all distinct, $r_{2,l'} \neq r_{1, l}$ for any $1 \leq l' \leq L'$ and $1 \leq l \leq L$,  $\beta: X_{r_{2,l}} \to X_{r_{2,l+1}}$, for $1 \leq l \leq L'-1$ and $\beta  : X_{r_{2,L'}} \to X_{r_{2,1}}$. 

Without loss of generality, we may assume that $L \leq L'$. Let $1 \leq c \leq N$ satisfy $\varphi^c : X_{r_{1,i}} \to X_{r_{2,2}}$. Then we have a diagram 

\begin{displaymath}
\xymatrix{ X_{r_{1,2}} \ar[r]^{\beta} \ar[d]^{\varphi^c} & X_{r_{1,2}} \ar[r]^{\beta}  \ar[d]^{\varphi^c} & \dots \ar[r]^{\beta}  & X_{r_{1,L}} \ar[r]^{\beta} \ar[d]^{\varphi^c} & X_{r_{1,1 \,}} \ar[d]^{\varphi^c}\\
X_{r_{2,2}} \ar[r]^{\beta} & X_{r_{2,2}} \ar[r]^{\beta}  & \dots \ar[r]^{\beta}& X_{r_{2,L}} \ar[r]^{\beta} & X_{r_{2,L+1}} \\
}
\end{displaymath}
which commutes. If follows that $X_{r_{2, L+1}} = Y_{r_{2, 1}}$, that is, $L = L'$. The proof now follows from induction.
\end{nproof}
\en

\bn
\begin{ntheorem} \label{IrrRok}
Suppose $(X, \varphi)$ is an irreducible Smale space with an effective $\mathbb{Z}$-action. Then the induced action on $\halg$ has finite Rokhlin dimension.
\end{ntheorem}

\begin{nproof}
Let $\beta : X \to X$ be the homeomorphism generating the $\mathbb{Z}$ action and let $X = X_1 \sqcup \cdots \sqcup X_N$ be the Smale decomposition into mixing components $(X_i, \varphi^N|_{X_i})$. Denote by $\mathrm{C^*}(H_i)$ the homoclinic algebra for the mixing Smale space $(X_i, \varphi^N|_{X_i})$. By Lemmas~\ref{sigma} and \ref{LS}, there exists $L$ such that 
\[ \beta^L  : X_i \to X_i \]
for every $i \in \{1, \dots, N\}$.  

For each $i$, we consider the action $\beta^L|_{X_i}$. Since $\beta$ induces an effective $\Z$-action on $(X, \varphi)$ (in particular $\beta^L \circ \varphi^N= \varphi^N\circ\beta^L$) it follows that $\beta^L|_{X_i}$ induces an effective $\Z$-action on the mixing Smale space $(X_i, \varphi^N|_{X_i})$. Hence, we can apply Theorem \ref{ZdActRok} to $\beta^L|_{X_i}$ acting on $(X_i, \varphi^N|_{X_i})$ to conclude that the action induced by $\beta^L|_{X_i}$ on $\mathrm{C}^*(H_i)$ has Rokhlin dimension $d_i$ for some $d_i < \infty$.

Let $\epsilon > 0$, $p \in \mathbb{N} \setminus \{0\}$. 

Let $f^{(k)}_{i, j}$, $i\in \{0,1\}$, $0 \leq j \leq p-1+i$, $0 \leq k \leq d_1$ be two Rokhlin towers for $\beta^L|_{X_1}$ with respect to $\epsilon$ and $p$ and any finite subset $\mathcal{F}_1  \subset \mathrm{C}^*(H_1)$.

We claim that the elements $\beta^l(f^{(k)}_{i, j})$, $1 \leq l \leq L-1$ satisfy (i) -- (v) of Lemma~\ref{lastConIsClearZ} with respect to $\halg, \epsilon$ and $p$.

Indeed, (i), (iii) and (iv) are clear by construction. To see (ii), we note that $\| \sum_{i,j,k} f^{(k)}_{i, j} - 1_{\mathrm{C}^*(H_1)}\| < \epsilon$ and since $\beta$ is a homeomorphism, it follows that $\| \sum_{i,j,k} \beta^l(f^{(k)}_{i, j}) - 1_{\mathrm{C}^*(H_l)})\| < \epsilon$ for every $l$, and hence  \[
\| \sum_{l,i,j,k} \beta^l(f^{(k)}_{i, j}) - 1_{\halg}\| < \epsilon.\]
\end{nproof}
\en

\bn
\begin{nremark}
Let $Y$ be a compact Hausdorff space. It is not difficult to check that if a $\mathbb{Z}$-action on $C(Y)$ has finite Rokhlin dimension, then the action on $Y$ must be free. Thus when $\Z$ acts on an irreducible Smale space $(X, \varphi)$, the action of $\Z$ on $C(X)$ cannot have finite Rokhlin dimension because the action of $\Z$ on $X$ is not free by Proposition \ref{notManyFreeActions}. Nevertheless, Proposition \ref{ZdActRok} implies that a $\Z$-action on $\halg$ induced from an effective action on $(X, \varphi)$ has Rokhlin dimension at most one.
\end{nremark}
\en

\section{The induced action on the stable and unstable algebras} \label{IndActSandUSec}

In this section we use what we have already proved about actions on the homoclinic algebra to deduce results for actions on the stable and unstable algebras. To do so, we take advantage of the embedding of the homoclinic algebra into the multiplier algebras of $\salg$ and $\ualg$. It is easy to check that if $(X, \varphi)$ is a Smale space, then $(X, \varphi^{-1})$ is also a Smale space with the bracket reversed. The unstable relation of $(X, \varphi)$ is then the stable relation of $(X, \varphi^{-1})$. Thus, it is enough to show something holds for $\salg$ of an arbitrary (irreducible, mixing) Smale space to imply the same for $\ualg$ of an arbitrary (irreducible, mixing) Smale space. Throughout this section, we again assume that $(X, \varphi)$ is irreducible, but is not necessarily mixing unless explicitly stated. 

\bn
Let $(X, \varphi)$ be an irreducible Smale space. Fix a finite set $P$ of $\varphi$-invariant periodic points. Let $\halg$ denote the associated homoclinic and $\salg$ the stable algebra. As in \cite[Theorem 3.4]{Put:C*Smale}, we define, for each $a \in C_c(\mathcal{G}_H)$, an element $(\rho(a), \rho(a))$ in $\mathcal{M}(\salg)$ by
\[ (\rho(a) b)(x,y) = \sum_{z \in X^U(P), \ z \sim_s x} a(x,z)b(z,y) \]
and
\[ (b \rho(a)) (x,y) = \sum_{z \in X^U(P), \ z \sim_s x} b(x,z) a(z,y)\]
where $b\in C_c(S)$. This extends to a map 
\[\rho : \halg \into \mathcal{M}(\salg).\]
We note that $\rho$ and the representations of these algebras on the Hilbert space $l^2(X^H(P,Q))$ (see Remark \ref{remRepSU}) are compatible. 
\en

\bn
\begin{nlemma} \label{MakeCommute}
Let $F \subset {C}_c(S)$ be a finite subset and let $r_1, \dots, r_N \in {C}_c(\mathcal{G}_H)$. Then, for every $\epsilon > 0$ there exists a $k\in \mathbb{N}$ such that, viewing $a$ as an element of the multiplier algebra $\mathcal{M}(\salg)$, we have
\[ \| \rho(\alpha^k(r_i)  ) a - a \rho(\alpha^k(r_i )) \| < \epsilon\]
for every $i = 1, \dots, N$ and every $a \in F$. 
\end{nlemma}

\begin{nproof}
Follows from \cite[Theorem 3.5]{Put:C*Smale} or \cite[Proposition 4.1]{Kil:thesis}.
\end{nproof}
\en

\bn
\begin{ndefn} 
Suppose $\beta : G \to \Aut(A)$ is an action of a countable discrete group. Let $I$ be a separable $G$-invariant ideal in $A$ and $B$ be a $\sigma$-unital $G$-$\mathrm{C}^*$-subalgebra of $A$. Then there exists a countable $G$-quasi-invariant quasicentral approximate unit $(w_n)_{n\in N}$ of $I$ in $B$. That is, there exists $(w_n)_{n\in \N}$ an approximate identity for $I$ such that 
\begin{enumerate}
\item for any $a \in B$, $\| aw_n - w_n a\| \rightarrow 0$ as $n\rightarrow \infty$;
\item for each $g\in G$, $\| \beta_g(w_n) -w_n \| \rightarrow 0$ as $n\rightarrow \infty$. 
\end{enumerate}
\end{ndefn}
\en

The existence of a $G$-quasi-invariant quasicentral approximate unit is shown in \cite{Kas:Novikov}, but also see \cite[5.3]{GHT:EthCstar} and \cite{Dum:BivTheCstar}.

\subsection{Rokhlin dimension for actions on the stable and unstable algebra}

To take advantage of the embedding of the homoclinic algebra into the multiplier algebras of $\salg$ and $\ualg$, we will find it convenient to introduce the  definition of  \emph{multiplier Rokhlin dimension with repect to a finite index subgroup}. It is inspired by \cite[Defintion 2.15]{Phi:Free}, \cite[Definition 1.1]{HirWinZac:RokDim} and \cite[Definitions 5.4, 10.2]{SzaWuZac:RP}. For a given $\mathrm{C}^*$-algebra $A$ and group $G$ with action $\beta : G \to {\rm Aut}(A)$ we denote by $\mathcal{M}(\beta)$ the induced action of $G$ on the multiplier algebra $\mathcal{M}(A)$. 

\bn
\begin{ndefn}\!\label{multiRD}
Let $A$ be a $\mathrm{C}^*$-algebra, $G$ a countable group with finite index subgroup $K$, $\beta : G \to \Aut(A)$ an action of $G$ on $A$ and $d \in \mathbb{N}$. We say the action has multilplier Rokhlin dimension $d$  with respect to $K$, denoted $\dimrok(\beta, K)$,  if $d$ is the least integer such that the following holds: for each $\epsilon>0$ and finite subsets $M \subset G$, $F \subset A$ there are positive contractions
$$ \left( f^{(l)}_{\overline{g}} \right)_{l=0, \ldots d; \overline{g}\in G/K} \subset \mathcal{M}(A) $$
such that
\begin{enumerate}
\item for any $l$, $\| f^{(l)}_{\overline{g}} f^{(l)}_{\overline{h}} \| < \epsilon$, for $\overline{g} \neq \overline{h}$ in $G/K$;
\item $\| \sum_{l=0}^d \sum_{\overline{g} \in G/K} f^{(l)}_{\overline{g}} - 1 \| < \epsilon$;
\item $\| \mathcal{M}(\beta_h)(f^{(l)}_{\overline{g}}) - f^{(l)}_{\overline{hg}} \| < \epsilon$  for all $l \in \{ 0, \ldots , d\}$, $\overline{g}  \in G/K$ and $h \in M$;
\item $\| [ f^{( l )}_{\overline{g}}, a] \| < \epsilon$ for all $l \in \{ 0, \ldots , d\}$, $\overline{g} \in G/K$ and $a \in F$.
\end{enumerate}
If no such $d$ exists, then we write $\dimrok(\beta, K) = \infty$. 
We say the action has multilplier Rokhlin dimension $d$ with commuting towers with respect to $K$, denoted $\dimrok^c(\beta, K)$,  if, in addition,
\begin{enumerate}[resume]
\item $\| f^{(l)}_{\overline{g}}f^{(k)}_{\overline{h}}- f^{(k)}_{\overline{h}}f^{(l)}_{\overline{g}} \| < \epsilon$ for every $k,l = 0, \dots, d$,  $\overline{g}, \overline{h} \in G/K$.
\end{enumerate}
If no such $d$ exists, then we write $\dimrok^c(\beta, K) = \infty$. 
\end{ndefn}
\en

The multiplier Rokhlin dimension with respect to a finite subgroup is used to define the Rokhlin dimension for an action of a countable residually finite group.

\bn
\begin{ndefn}[cf. \cite[Definition 5.8]{SzaWuZac:RP}] 
Let $A$ be a $\mathrm{C}^*$-algebra, $G$ a countable, residually finite group, and $\beta : G \to \Aut(A)$ an action of $G$ on $A$. The Rokhlin dimension of $\beta$ is defined by
\[ \dimrok(\beta) := \sup \{\dimrok(\beta, K) \mid K \leq G, [G:H] < \infty\}.\]
The Rokhlin dimension of $\beta$ with commuting towers is given by
\[ \dimrok^c(\beta) := \sup \{\dimrok^c(\beta, K) \mid K \leq G, [G:H] < \infty\}.\]
\end{ndefn}
\en

It will be easier to show finite multiplier Rokhlin dimension relative to a finite subgroup, which we show implies the definition of Rokhlin dimension relative to a finite index subgroup for $\mathrm{C}^*$-algebras given in \cite{SzaWuZac:RP}, recalled in Definition~\ref{RDwrtH} below. Some further notation is required to do so.  

First, we need to define central sequence algebras. Loosely speaking, working in a central sequence algebra allows one to turn statements such as approximate commutativity in the original algebra into honest commutativity in the central sequence algebra. As such, central sequence arguments often allow one to streamline proofs. In the case of discrete groups, an action on a $\mathrm{C}^*$-algebra induces an action on its central sequence algebra, and it is possible to reformulate definitions for both Rokhlin dimension of finite group actions and integer actions on separable unital $\mathrm{C}^*$-algebras in terms of induced actions on central sequence algebras \cite[Lemma 2.4]{SzaWuZac:RP}. 

\bn
\begin{ndefn}
 Let $A$ be a separable $\mathrm{C}^*$-algebra. We denote the sequence algebra of $A$ by
\[ A_{\infty} := \prod_{n \in \mathbb{N}} A / \bigoplus_{n \in \mathbb{N}} A. \]
We view $A$ as a subalgebra of $A_{\infty}$ by mapping an element $a \in A$ to the constant sequence consisting of $a$ in every entry. The central sequence algebra is then defined to be
\[ A^{\infty} := A_{\infty} \cap A' = \{x \in A_{\infty} \mid ax = xa \text{ for every } a \in A\}, \]
the relative commutant of $A$ in $A_{\infty}$. Let
\[ \mathrm{Ann}(A, A_{\infty}) := \{ x \in A_{\infty} \mid ax = xa = 0 \text{ for every } a \in A\},\]
which is evidently an ideal in $A^{\infty}$. Finally, we define
\[ F(A) := A^{\infty}/ \mathrm{Ann}(A, A_{\infty}). \]
\end{ndefn}
\en

When $A$ is not separable, one can define the above with respect to a given separable subalgebra $D$, as is done in \cite{SzaWuZac:RP}. However, since all our $\mathrm{C}^*$-algebras will be separable, we will not require this. 

A completely positive contractive (c.p.c.) map $\varphi : A \to  B$ between $\mathrm{C}^*$-algebras $A$ and $B$ is said to be \emph{order zero} if it is orthogonality preserving, that is, for every $a, b \in A_+$ with $ab = ba =0$ we have $\varphi(a) \varphi(b) = 0$. Any $^*$-homomorphism is of course order zero, but a c.p.c. order zero map is not in general a $^*$-homomorphism. For more about c.p.c. order zero maps, see \cite{WinZac:order-zero}.

\bn
\begin{ndefn}\label{RDwrtH}\cite[Definitions 5.4, 10.2]{SzaWuZac:RP}
Let $A$ be a separable $\mathrm{C}^*$-algebra, $G$ a countable group with finite index subgroup $K$, $\beta : G \to \Aut(A)$ an action of $G$ on $A$ and $d \in \mathbb{N}$. Let $\tilde{a}_{\infty}$ denote the induced action on $F_{\infty}(A)$. We say the action has Rokhlin dimension $d$  with respect to $K$ if $d$ is the least integer such that there exists equivariant c.p.c. order zero maps 
\[ \varphi_l : (C(G/K), G\text{-shift}) \to (F_{\infty}(A), \tilde{a}_{\infty}), \quad l = 0, \dots, d\]
with 
\[ \varphi_0(1) + \cdots + \varphi_d(1) = 1. \]
If moreover $\varphi_0 ,\dots, \varphi_d$ can be chosen to have commuting ranges, then we say the action has Rokhlin dimension $d$  with commuting towers with respect to $K$.
\end{ndefn}
\en

The next result is an obvious generalisation of the equivalence of (1) and (3) of \cite[Proposition 4.5]{SzaWuZac:RP} to the case of commuting towers. 

\bn
\begin{nlemma} \label{ComTow} Let $A$ be a $\mathrm{C}^*$-algebra, $G$ a countable group and $K$ a subgroup of finite index. The following are equivalent for an action $\beta : G \to \Aut(A)$.
\begin{enumerate}
\item The action $\beta$ has Rokhlin dimension $d$ with commuting towers with respect to $K$.
\item For every finite subset $M \subset G$, finite subset $F \subset A$ and $\epsilon >0$ there are positive contractions $(f^{(l)}_{\overline{g}})_{\overline{g} \subset H}^{l= 0, \dots, d}$ in $A$ such that
	\begin{enumerate}
	\item $\|(\sum_{l =0}^d \sum_{\overline{g} \in G/K} f^{(l)}_{\overline{g}}) \cdot a - a \| < \epsilon$ for all $a \in F$;
	\item $\| f_{\overline{g}}^{(l)}f_{\overline{h}}^{(l)}a \| \leq \epsilon$ for all $a\in F$, $l\in \{0, \dots, d\}$ and $\overline{g} \neq \overline{h} \in G/K$;
	\item $\|(\beta_g(f^{(l)}_{\overline{h}}) - f^{(l)}_{\overline{gh}})a \| < \epsilon$ for all $a \in F$, $l \in 0, \dots, d$ and $g \in M$ and $\overline{h} \in G/K$;
	\item $\| f^{(l)}_{\overline{g}}a - a f^{(l)}_{\overline{g}} \| < \epsilon$ for all $a \in F$, $l\in \{0, \dots, d\}$ and $\overline{g} \in G/K$;
	\item $\| (f^{(k)}_{\overline{g}} f^{(l)}_{\overline{h}} - f^{(l)}_{\overline{h}} f^{(k)}_{\overline{g}})a \| < \epsilon$ for all $a \in F$, $k,l\in \{0, \dots, d\}$ and $\overline{g}, \overline{h} \in G/K$.
	\end{enumerate}
\end{enumerate}
\end{nlemma} 
\begin{nproof} The only thing that one needs to check is that asking for the $f^{(l)}_{\overline{g}}$ to approximately commute is equivalent to having the images of the order zero maps of \cite[Definition 4.4]{SzaWuZac:RP} commute, but this is obvious.
\end{nproof}
\en

\bn
\begin{ntheorem} \label{MultToRok}
Let $G$ be a countable discrete group, $K$ a subgroup of $G$ with finite index, $A$ a separable $\mathrm{C}^*$-algebra and $\beta : G \to \Aut(A)$ an action with multiplier Rokhlin dimension at most $d$ with respect to $K$. Then $\beta$ has Rokhlin dimension at most $d$ respect to $K$.

If $\beta$ has multiplier Rokhlin dimension at most $d$ with commuting towers with respect to $K$, then $\beta$ has Rokhlin dimension at most $d$ with commuting towers respect to $K$.
\end{ntheorem}

\begin{nproof} We show the action satisfies the criteria of \cite[Proposition 5.5 (3)]{SzaWuZac:RP} (and in the commuting tower case Lemma~\ref{ComTow} (ii)). Note that \cite[Proposition 5.5 (3)]{SzaWuZac:RP} is exactly (a)-(d) in Lemma~\ref{ComTow} (ii).

Let $M \subset G$ and $F \subset A$ be finite subsets and let $\epsilon > 0$. Without loss of generality we may assume that every $a \in F$ is a positive contraction. Since $\beta$ has multiplier Rokhlin dimension less than or equal to $d$ with respect to $K$ we can find positive contractions $(f_g^{(l)})_{\overline{g} \in G/K}$, $l = 0 , \dots d$ satisfying Definition~\ref{multiRD}  with respect to $M$, $F$ and $\epsilon/2$. 

Let $(w_n)_{n \in \mathbb{N}}$ be an $G$-quasi-invariant quasicentral approximate unit for $A$ in $\mathcal{M}(A)$. Since $F \subset A$ and $(f_g^{(l)})_{\overline{g} \in G/K}$, $l = 0 , \dots d$ are in $\mathcal{M}(A)$, there exists $N \in \mathbb{N}$ sufficiently large so that
\[ \| w_N a - a \| < \epsilon/4 \text{ for every } a \in F,\]
and 
\[ \| [f_{\overline{g}}^{(l)}, w_N] \| < \epsilon/8.\]

and
\[ \| \beta_g(f_{\overline{h}}^{(l)}) - f_{\overline{gh}}^{(l)} \| < \epsilon/8.\]

By increasing $N$ if necessary, we may also assume, since $M$ is finite, that
\[ 
\| \beta_g(w_N) - w_N \| < \epsilon/8
\]
 for each $g\in M$.

Then let
\[ r_{\overline{g}}^{(l)} := w_N f_{\overline{g}}^{(l)} w_N,\]
for $l \in \{ 0, \dots, d\}$ and $g \in G$. Then each $r_{\overline{g}}^{(l)}$ is a positive contraction in $A$ and 
\begin{eqnarray*}
\| (\tsum_{l=0}^d\tsum_{\overline{g} \in G/K} r^{(l)}_{\overline{g}}) a - a\| &=& \| (\tsum_{l=0}^d\tsum_{\overline{g} \in G/K} w_N f^{(l)}_{\overline{g}} w_N) a - a\|  \\
&\leq& \| \tsum_{l=0}^d\tsum_{\overline{g} \in G/K} w_N f^{(l)}_{\overline{g}} w_N - 1 \|\\
&\leq&  \| \tsum_{l=0}^d\tsum_{\overline{g} \in G/K} f^{(l)}_{\overline{g}} - 1 \|\\
&<& \epsilon,
\end{eqnarray*}
showing that (a) of  Lemma~\ref{ComTow} (ii) holds. 

Next, 
\begin{eqnarray*}
\| r_{\overline{g}}^{(l)} r_{\overline{h}}^{(l)} a \| &=& \| w_N f_{\overline{h}}^{(l)} w_N^2 f_{\overline{g}}^{(l)} w_N a \| \\
&=& \epsilon/2 +  \| w_N^2 f_{\overline{h}}^{(l)} f_{\overline{g}}^{(l)} w_N^2 a \| \\
&<& \epsilon
\end{eqnarray*}
for every $a \in F$, $l\in \{0, \dots, d\}$ and $\overline{g} \in G/K$, showing that (b) of  Lemma~\ref{ComTow} holds.

Using $\| \beta_g(w_N) - w_N \| < \epsilon/8$ for each $g\in M$, we obtain
\begin{eqnarray*}
 \| (\beta_g(r_{\overline{h}}^{(l)}) - r_{\overline{gh}}^{(l)})a\| &=&  \epsilon/4 + \| w_N \mathcal{M}(\beta_g)(f^{(l)}_{\overline{h}}) w_N - w_N f^{(l)}_{\overline{gh}} w_N \| \\
 &<& \epsilon,
 \end{eqnarray*}
 for every $l \in \{0, \dots, d\}$, every $\overline{h} \in G/K$, every $g \in M$ and every $a \in F$, showing (c) of  Lemma~\ref{ComTow}.

For (d) of Lemma~\ref{ComTow} we have
\[
\| r_{\overline{g}}^{(l)} a - a r_{\overline{g}}^{(l)} \| = \epsilon/ 2 + \| f_{\overline{g}}^{(l)} a - a f_{\overline{g}}^{(l)} \| < \epsilon,
\]
for every $a \in F$, every $l \in \{0, \dots, d\}$ and every $\overline{g} \in G/K$.

Finally, in the commuting tower case, we must show (e): 
\begin{eqnarray*}
\| (r_{\overline{g}}^{(k)} r_{\overline{h}}^{(l)} -  r_{\overline{h}}^{(l)} r_{\overline{g}}^{(k))} a\|  &\leq& \| f_{\overline{g}}^{(k)} w_N^2 f_{\overline{h}}^{(l)} -  f_{\overline{h}}^{(l)} w_N^2  f_{\overline{g}}^{(k)}\| \\
&\leq& \epsilon/4 + \| f_{\overline{g}}^{(k)}  f_{\overline{h}}^{(l)} -  f_{\overline{h}}^{(l)}  f_{\overline{g}}^{(k)}\|\\
&<& \epsilon,
\end{eqnarray*}
for every $\overline{g}, \overline{h} \in G/K$ and $k, l \in \{0, \dots, d\}$.
\end{nproof}
\en

The proof of the next lemma is obvious and hence omitted. It will, however, prove useful in what follows.

\bn
\begin{nlemma} Suppose that the action of $G$ on $\halg$ has Rokhlin dimension at most $d$. Then we may choose the Rokhlin elements to satisfy $f^{(l)}_g \in C_c(\mathcal{G}_H)$ for  $l = 0, \dots, d$ and $ g \in G$.
\end{nlemma}
\en

\bn
\begin{nprop} \label{G-RD for S}
Let $G$ be a countable group acting on an irreducible Smale space $(X, \varphi)$ and let $K \subset G$ be a subset of finite index. Then if the induced action $\beta : G \to \Aut(\halg)$ has Rokhlin dimension at most $d$ with respect to $K$ so does the action $\beta^{(S)} : G \to \Aut(\salg)$. 

If $\beta$ has Rokhlin dimension at most $d$ with commuting towers with respect to $K$ so does the action $\beta^{(S)}$.
\end{nprop}

\begin{nproof} We will show that $\beta^{(S)}$ has multiplier Rokhlin dimension at most $d$ with respect to $K$. The result then follows from Theorem~\ref{MultToRok}. 

Let $M$ be a finite subset of $G$, $F$ a finite subset of $\salg$ and $\epsilon >0$. Without loss of generality, we may assume that $F \subset C_c(S)$.  Since $\beta$ has Rokhlin dimension at most $d$ with respect to $K$ there are contractions $r^{(l)}_g \in \halg$ such that 
\begin{enumerate}
\item $\| r^{(l)}_{\overline{g}} r^{(l)}_{\overline{h}}  \| < \epsilon/2$ for $l = 0, \dots, d$ and any $g, h \in G$ with $g \neq h$,
\item $\| \tsum_{l=0}^d \tsum_{\overline{g} \in G/K} r^{(l)}_{\overline{g}} - 1\| < \epsilon/2,$
\item $\| \beta_h(r^{(l)}_{\overline{g}}) - r^{(l)}_{\overline{hg}}\| < \epsilon/2$ for $l =0, \dots, d$, every $h \in M$ and $\overline{g} \in G/K$, and
\item $\| r^{(k)}_{\overline{g}}r^{(l)}_{\overline{h}} - r^{(l)}_{\overline{h}}r^{(k)}_{\overline{g}}\| < \epsilon/2$ for every $\overline{g}$, $\overline{h} \in G/K$ and $k,l \in \{0, \dots, d\}$.
\end{enumerate}
By the previous lemma, we may moreover  assume that each $r^{(l)}_{\overline{g}}  \in C_c(\mathcal{G}_H)$. Now we can find a natural number $k$ such that 
\[ \| \rho(\alpha^k(r^{(l)}_{\overline{g}})) a - a \rho(\alpha^k(r^{(l)}_{\overline{g}})) \| <  \epsilon/2,\]
for every $a \in F$. 

For $l \in \{0, \dots, d\}$ and $\overline{g} \in G/K$, let 
\[ f^{(l)}_{\overline{g}} := \rho(\alpha^k(r^{(l)}_{\overline{g}} )).\]
Note that each $f^{(l)}_g$ is a positive contraction in $\mathcal{M}(\salg)$. We will show that the $f^{(l)}_{\overline{g}}$  satisfy (i) -- (iv) of Definition~\ref{multiRD}. 

For (i) we have 
\begin{eqnarray*}
 \| f^{(l)}_{\overline{g}} f^{(l)}_{\overline{h}} \| &=&  \|  \rho(\alpha^k(r^{(l)}_{\overline{g}})) \rho(\alpha^k(r^{(l)}_{\overline{h}} )) \| \\
  &=&  \| \alpha^k(r^{(l)}_{\overline{g}} r^{(l)}_{\overline{h}}) \| \\
  &<& \epsilon,
 \end{eqnarray*}
 for any $\overline{g} \neq \overline{h} \in G/K$, any $l \in \{0, \dots, d\}$.
 
For (ii)
\begin{eqnarray*}
 \| \tsum_{l=0}^d \tsum_{\overline{g} \in G/K} f_{\overline{g}}^{(l)} - 1_{\mathcal{M}(\salg)} \| &=& \|  \tsum_{l=0}^d \tsum_{\overline{g}\in G} \rho(\alpha^k( r^{(l)}_{\overline{g}} )) - \rho(1_{\halg})\| \\
 &=& \| \tsum_{l=0}^d \tsum_{\overline{g} \in G/K} r^{(l)}_{\overline{g}} - 1_{\halg} \| \\
 &<& \epsilon, 
 \end{eqnarray*}
 for any $a \in F$.
 
 For (iii) we have
 \begin{eqnarray*} 
 \| f^{(l)}_{\overline{g}} a - a f^{(l)}_{\overline{g}} \| &=& \| \rho(\alpha^k(r^{(l)}_{\overline{g}})) a - a \rho(\alpha^k(r^{(l)}_{\overline{g}})) \| \\
 &<&  \epsilon/2,
 \end{eqnarray*}
 for all $l \in \{0, \dots d\}, \overline{g} \in G/K$ and $a \in F$.

For (iv), let $g \in G$ with $\overline{h} \in G/K $ and let $a \in F$. Then, 
\begin{eqnarray*}
\| \beta^{(S)}_g (f^{(l)}_{\overline{h}}) - f^{(l)}_{\overline{gh}} \| &=& \| \alpha^k(\beta_g(r^{(l)}_{\overline{h}}) -  r^{(l)}_{\overline{gh}})  \| \\
&<& \epsilon.
\end{eqnarray*}
Finally, in the commuting towers case, for (v), let $\overline{g}, \overline{h} \in G/K$ and $l,k \in \{0, \dots, d\}$,
\begin{eqnarray*}
\| f^{(k)}_{\overline{h}} f^{(l)}_{\overline{g}} -  f^{(l)}_{\overline{g}}f^{(k)}_{\overline{h}}  \| &=& \| \rho(\alpha^k(r^{(k)}_{\overline{g}}  r^{(l)}_{\overline{h}})  - \rho(\alpha^k(r^{(l)}_{\overline{h}} r^{(k)}_{\overline{g}}) \| \\
&=& \| r^{(k)}_{\overline{g}}  r^{(l)}_{\overline{h}} - r^{(l)}_{\overline{h}}  r^{(k)}_{\overline{g}}\|  \\
&<& \epsilon. 
\end{eqnarray*}
The result now follows.
\end{nproof}
\en

This gives us the next corollary:

\bn
\begin{ncor} Let $G$ be a countable residually finite group acting on an irreducible Smale space $(X, \varphi)$. Then if the induced action $\beta : G \to \Aut(\halg)$ has Rokhlin dimension at most $d$ (with commuting towers) so does the action $\beta^{(S)} : G \to \Aut(\salg)$.
\end{ncor}
\en

By Proposition \ref{ZdActRok} and Theorem~\ref{IrrRok} respectively, we get the next two corollaries.

\bn
\begin{ncor}
Let $d\in \N$ and $G=\Z^d$. Suppose $G$ acts effectively on a mixing Smale space $(X, \varphi)$. Then the induced actions $\beta : G \to \Aut(\halg)$, $\beta^{(S)} : G \to \Aut(\salg)$ and $\beta^{(U)} :  G \to \Aut(\ualg)$ each have finite Rokhlin dimension.
\end{ncor}
\en

\bn
\begin{ncor} Let $\mathbb{Z}$ be an effective action on an irreducible Smale space $(X, \varphi)$. Then the induced actions $\beta : \mathbb{Z} \to \Aut(\halg)$, $\beta^{(S)} :  \mathbb{Z} \to \Aut(\salg)$ and $\beta^{(U)} :  \mathbb{Z} \to \Aut(\ualg)$ each have finite Rokhlin dimension.
\end{ncor}
\en

For more general group actions the situation is less clear. In particular, to the authors' knowledge, it is not known if strong outerness implies finite Rokhlin dimension\footnote{At the workshop ``Future Targets in the Classification Program for Amenable $\mathrm{C}^*$-Algebras'' held at BIRS, Eusebio Gardella presented results in this direction, but they have yet to appear.}.

In fact, for general discrete groups, there are obstructions to a (strongly outer) action having finite Rokhlin dimension with commuting towers, see \cite[Corollary 4.6]{HirPhi:RD}. Such examples can occur in the context considered in the present paper. An explicit example is the following, let $(\Sigma_{[3]}, \sigma)$ be the full three shift (so $\Sigma_{[3]}= \{ 0, 1, 2\}^{\Z}$ and $\sigma$ is the left sided shift). Then $\halg$ is the UHF-algebra with supernatural number $3^{\infty}$. The action induced from the permutation $0 \mapsto 1$, $1 \mapsto 0$, and $2 \mapsto 2$ is an effective order two automorphism. Thus the action induced on $\halg$ is strongly outer, but it follows from \cite[Corollary 4.6 (2)]{HirPhi:RD} that it does not have finite Rokhlin dimension with commuting towers.

\section{$\mathcal{Z}$-stability, nuclear dimension and classification}

We begin this section with two theorems that follow quickly from the work done above. As well as being interesting observations on their own, they will  allow us to say something about the $\mathcal{Z}$-stability of the crossed products of $\salg$ and $\ualg$ below.

\bn
\begin{ntheorem} \label{HZstab}
Let $G$ be a countable discrete amenable group. Suppose $G$ acts on a mixing Smale space $(X, \varphi)$. Denote by $\beta : G \to \Aut(\halg)$ the induced action of $G$ on $\halg$. Then $\halg \rtimes_{\beta} G$ is $\mathcal{Z}$-stable. 
\end{ntheorem}

\begin{nproof} Since $(X, \varphi)$ is mixing, $\halg$ is simple and so the classification results of \cite{DS:NucDimSmale} imply that $\halg$ is $\mathcal{Z}$-stable.  Since $\halg$ has unique trace,  it must be fixed by the action. Thus $\halg \rtimes_{\beta} G$ is $\mathcal{Z}$-stable by \cite[Theorem 1.1]{Sat:Amen}.
\end{nproof}
\en

\bn
\begin{ntheorem} Let $G$ be a countable amenable group acting effectively on a mixing Smale space $(X, \varphi)$ and $\beta : G \to \Aut(\halg)$ the induced action. Then the crossed product $\halg \rtimes_{\beta} G$ is a simple unital nuclear $\mathcal{Z}$-stable $\mathrm{C}^*$-algebra with unique tracial state and nuclear dimension (in fact, decomposition rank) at most one.  In particular, it belongs to the class of $\mathrm{C}^*$-algebras that are classified by the Elliott invariant and $(\halg \rtimes_{\beta} G) \otimes \uhf$ is TAF, for any UHF algebra of infinite type $\uhf$ .
\end{ntheorem}

\begin{nproof}
Since $\halg$ is amenable, it follows from \cite[II.3.9]{Ren:groupoid} that $\halg \rtimes_{\beta} G \cong \mathrm{C}^*(\mathcal{G}_H \rtimes G)$ is also amenable and hence by \cite{Tu:Groupoids} it satisfies the UCT. Thus  $\halg \rtimes_{\beta} G$ is quasidiagonal  \cite[Corollary B]{TikWhiWin:QD} and by \cite[Corollary E]{TikWhiWin:QD} together with Theorem~\ref{HZstab} classified by the Elliott invariant.  The nuclear dimension and decomposition rank bounds are given by \cite[Theorem F]{BBSTWW:2Col}.  Finally,  by \cite[Theorem 6.1]{MatSat:UHF}  $(\halg \rtimes_{\beta} G) \otimes \uhf$ is TAF. 
\end{nproof}
\en

\subsection{$\mathcal{Z}$-stability of crossed products $\salg$ and $\ualg$}

Let $A$ be separable $\mathrm{C}^*$-algebra and $G$ a discrete group.

For an action $\beta: G\to \Aut(A)$ the fixed point algebra is given by

\[ A^{\beta} := \{ a \in A \mid \beta_g(a) = a \text{ for all } g \in G\}. \]

Let
\[ A^{\infty} := \ell^{\infty}(\mathbb{N}, A) / c_0(A) .\]

The central sequence algebra of $A$ is defined by
\[ A_{\infty} := A^{\infty} \cap A',\]
where $A$ is considered as the subalgebra of $A^{\infty}$ by viewing an element as a constant sequence.

We will denote by $\overline{\beta}$ the induced action on $A_{\infty}$.

Let $p$ and $q$ be positive integers. The dimension drop algebra $I(p,q)$ is given by
\[
I(p, q) := \{ f \in C([0,1], M_p(\C) \otimes M_q(\C) \mid f(0) \in M_p(\C) \otimes \C \hbox{ and }f(1)\in \C \otimes M_q(\C)\}.
\]
The Jiang--Su algebra, $\mathcal{Z}$, is an inductive limit of such algebras \cite{JiaSu:Z}.

\bn
\begin{ntheorem} \label{zstab} Suppose $G$ is a countable discrete group acting  on a mixing Smale space $(X, \varphi)$. Let $\beta^{(S)}$ and $\beta^{(U)}$ denote the induced actions on $\salg$ and $\ualg$, respectively. If, for each $k\in \N$, there exists a unital equivariant embedding $I(k, k+1) \rightarrow \halg$, then $\salg \rtimes_{\beta^{(S)}} G$ and $\ualg \rtimes_{\beta^{(U)}} G$ are $\mathcal{Z}$-stable.
\end{ntheorem}

\begin{nproof}
As usual, it suffices to prove the result for $\salg \rtimes_{\beta^{(S)}} G$. To do so, we will show that the hypotheses of Lemma 2.6 in \cite{HirPhi:RD} hold. We will show that, for every $k \in \mathbb{N}$, there is a completely positive contractive map 
\[ \gamma : I(k, k+1) \to \salg_{\infty} \]
satisfying 
\begin{enumerate}
\item $ (\overline{\beta}^S( \gamma(x)) - \gamma(x))a = $ for every $x \in I(k, k+1)$, $g \in G$ and $a \in \salg$,
\item $a \gamma(1) = a$ for every $a \in \salg$, and 
\item  $a(\gamma(xy) - \gamma(x)\gamma(y)) = 0$ for every $x, y \in I(k, k+1)$ and $a \in \salg$.  
\end{enumerate}

Suppose $F \subset I(k, j+1)$ a finite subset and $\epsilon > 0$ are given. Let $w_n$ be a $G$-invariant quasicentral approximate unit for $\salg$ in $\mathcal{M}(\salg)$ and 
\[ \tilde{\gamma} : I(k, k+1) \to (\halg_{\infty})^{\beta}. \]
be a unital embedding (which exists by assumption).

Define
\[ \gamma : I(k, k+1) \to \salg^{\infty} \]
via
\[ \gamma(d) = (w_n \rho(\alpha^n(d_n)) w_n)_{n \in \mathbb{N}}, \]
where $(d_n)_{n \in \mathbb{N}}$ is a representative sequence for $\tilde{\gamma}(d)$. Then $\gamma$ gives a c.p.c.\! map. Moreover, if $a \in \salg$ we have   
\[  \lim_{n \to \infty} \| w_n \rho(\alpha^n(d_n) )w_n a - a w_n \rho(\alpha^n(d_n)) w_n \| = 0,\]
so in fact
\[  \gamma : I(k, k+1) \to \salg_{\infty}. \]
Let us check that $\gamma$ satisfies (i), (ii) and (iii).

Let $a \in A$ and for any $d \in I(k, k+1)$, let $(d_n)_{n \in\mathbb{N}}$ be a representative of $\tilde{\gamma}(d)$ in  $(\halg_{\infty})^{\beta}$.
\begin{eqnarray*}
\|  \overline{\beta}^S(\gamma(d) - \gamma(d)) a \|  &=&  \lim_{n \to \infty} \| \beta^{(S)}(w_n \rho(\alpha^n(d_n)) w_n) - w_n \rho(\alpha^n(d_n)) w_n\| \\
&=& \lim_{n \to \infty}  \| w_n \rho(\beta(\alpha^n(d_n)) w_n - w_n \rho(\alpha^n(d_n)) w_n\| \\
&\leq&  \lim_{n \to \infty} \| \beta(\alpha^n(d_n) ) - \alpha^n(d_n) \|\\
&=& 0,
\end{eqnarray*}
showing $(i)$.

To show (ii), we have
\begin{eqnarray*}
a \gamma(1) &=& (a w_n \rho(\tilde{\gamma}(1)) w_n)_{n \in \mathbb{N}}\\
&=& (a \rho(1) w_n^2)_{n \in \mathbb{N}} \\
&=& a,
\end{eqnarray*}
 for every $a \in \salg$

Finally, for (iii), let $d, d' \in I(k, k+1)$ and $a \in \salg$. Then
\begin{eqnarray*}
a(\gamma(d d') - \gamma(d)\gamma(d')) &=& (a w_n \rho(d_n d'_n) w_n)_{n \in \mathbb{N}} - (a w_n \rho(d_n) w_n^2 \rho(d'_n) w_n)_{n \in \mathbb{N}}\\
&=& (a w_n^2 \rho(d_n d_n'))_{n \in \mathbb{N}} - (a w_n^4 \rho(d_n d_n'))_{n \in \mathbb{N}} \\
&=& (a   \rho(d_n d_n'))_{n \in \mathbb{N}} -  (a   \rho(d_n d_n'))_{n \in \mathbb{N}}\\
&=& 0.
\end{eqnarray*}
Thus $\salg \rtimes_{\beta^{(S)}} G$ is $\mathcal{Z}$-stable.
\end{nproof}
\en

\bn \begin{ncor} Suppose $G$ is a discrete amenable group and $G$ acts on a mixing Smale space $(X, \varphi)$. Let $\beta^{(S)}$ and $\beta^{(U)}$ denote the induced actions on $\salg$ and $\ualg$, respectively. Then $\salg \rtimes_{\beta^{(S)}} G$ and $\ualg \rtimes_{\beta^{(U)}} G$ are $\mathcal{Z}$-stable.
\end{ncor}

\begin{nproof}
By Theorem~\ref{knownprops} (v), $\halg$ is $\mathcal{Z}$-stable and hence has strict comparison \cite[Corollary 4.6]{Ror:Z-absorbing}. This in turn implies $\halg$ has property (TI) of Sato \cite[Proposition 4.4]{Sat:Amen}. Since $\halg$ has unique trace $\tau$ which is therefore fixed by $\beta$, the hypotheses of \cite[Corollary 3.6]{Sat:Amen} are satisfied. Thus from the proof of \cite[Theorem 5.2]{Sat:Amen} we get a unital embedding 
\[ \tilde{\gamma} : I(k, k+1) \to (\halg_{\infty})^{\beta}. \]
The result then follows from the previous theorem.
\end{nproof}
\en

\

\emph{Acknowledgments.} The authors thank Ian Putnam for many useful discussions concerning the content of this paper, Smale spaces, group actions and dynamics in general. We also thank Magnus Goffeng for a number of useful comments. The authors thank the referee for reading the paper carefully and making a number of useful comments. The authors wish to thank Matrix at the University of Melbourne for hosting them during the programme \emph{Refining $\mathrm{C}^*$-algebraic Invariants for Dynamics using KK-theory} in July 2016, the Banach Centre at the Institute of Mathematics of the Polish Academy of Sciences for hosting the first listed author during the conference \emph{Index Theory} in October 2016, the University of Hawaii, Manoa for hosting the second listed author during the workshop \emph{Computability of K-theory} in November 2016 and the Centre Rercerca Matem\`atica, Barcelona, for their stay during the \emph{Intensive Research on Operator Algebras: Dynamics and Interactions} in July 2017. The above research visits were partially supported through NSF grants DMS 1564281 and DMS 1665118.

\end{document}